\documentclass[lettersize,journal]{IEEEtran}
\usepackage{amsmath,amsfonts}
\usepackage{algorithmic}
\usepackage{algorithm}
\usepackage{array}
\usepackage[caption=false,font=normalsize,labelfont=sf,textfont=sf]{subfig}
\usepackage{textcomp}
\usepackage{stfloats}
\usepackage{url}
\usepackage{verbatim}
\usepackage{graphicx}
\usepackage{cite}
\usepackage{mathrsfs}
\usepackage{color}
\usepackage{bm}
\usepackage{amssymb}
\usepackage{amsthm}
\usepackage{caption}
\usepackage{epstopdf}
\usepackage{booktabs}
\usepackage{epsfig}
\usepackage{float}
\usepackage{xcolor}  

\usepackage{hyperref}
\usepackage{todonotes}

\begin{document}
		
\title{Two-Channel Filter Banks on Joint Time-Vertex Graphs with Oversampled Graph Laplacian Matrix}
		
\author{Yu Zhang and Bing-Zhao Li$^{\ast}$,~\IEEEmembership{Member,~IEEE}
\thanks{The current study was partially supported by the National Natural Science Foundation of China (No. 62171041), and Natural Science Foundation of Beijing Municipality (No. 4242011), \textit{(Corresponding author: Bing-Zhao Li.)}}
\thanks{Yu Zhang is with the School of Mathematics and Statistics, Beijing Institute of Technology, Beijing 102488, China (e-mail: zhangyu$\_$bit@outlook.com).}
\thanks{Bing-Zhao Li is with the School of Mathematics and Statistics, Beijing Institute of Technology, Beijing 102488, China (e-mail: li$\_$bingzhao@bit.edu.cn).}}

\markboth{Journal of \LaTeX\ Class Files,~Vol.~, No.~, ~2021}%
{Shell \MakeLowercase{\textit{et al.}}: A Sample Article Using IEEEtran.cls for IEEE Journals}
		
\IEEEpubid{0000--0000/00\$00.00~\copyright~2021 IEEE}
		
\maketitle
		
\begin{abstract}
To address the limitations of conventional critically sampled graph filter banks in joint time-vertex signal processing, which require decomposing the joint graph into bipartite subgraphs and thus cannot fully exploit all temporal and spatial edges in a single-stage transform, we introduce the joint time-vertex oversampled graph Laplacian matrix. This operator enables the construction of bipartite extensions that preserve all edges of the original joint graph and supports redundant multiresolution representations. Based on this operator, we design two-channel joint time-vertex oversampled graph filter banks and develop efficient oversampling extensions using a $K$-coloring strategy. The proposed framework is applied to both graph signal and image/video denoising, modeling images as graph signals to leverage structural relationships. Extensive experiments demonstrate its effectiveness in decomposition, reconstruction, and denoising, achieving notable performance improvements over critically sampled and existing methods.
\end{abstract}
\begin{IEEEkeywords}
Graph filter banks, joint time-vertex signal processing, graph oversampling, graph Laplacian matrix, denoising.
\end{IEEEkeywords}


\section{Introduction}
\label{Intro}
\IEEEPARstart{W}{ITH} the rapid advancement of information technology, the global volume of stored and processed data has been growing exponentially. A significant portion of these data are generated on complex networks with inherently irregular structures, which are naturally modeled as graphs, while the associated data are represented as time-varying graph signals. In recent years, graph signal processing (GSP) has evolved rapidly, establishing itself as a core framework for analyzing and processing such data \cite{GFTlaplace,GFTadjacency1,GFTadjacency2,Goverview,Gwavelet, Ghistory,Gvertex,Gfrequency,Gfilters,GFTsampling,GUncertainty,GFTuncertainty}.

Graphs serve as a powerful tool for representing irregular and complex data and have found widespread applications in social networks, sensor networks, transportation systems, as well as image and video processing \cite{Gengraphs,Aligraphs,multi-agentgraphs}. Within this paradigm, signals are mapped to graph vertices, enabling GSP to provide effective tools for network data processing. Considerable efforts have been devoted to extending classical time-frequency methods to the graph domain, including discrete Fourier transform (DFT) \cite{GFTlaplace,GFTadjacency1,GFTadjacency2,Gwavelet}, windowed Fourier transform \cite{Gvertex}, spectral analysis \cite{Gfrequency}, filtering \cite{Gfilters}, sampling and interpolation \cite{GFTsampling,GUncertainty}, and uncertainty principles \cite{GUncertainty,GFTuncertainty}. Analogous to the role of DFT \cite{DFT} in conventional signal processing, the graph Fourier transform (GFT) has become fundamental to many GSP applications, such as denoising \cite{biorth,PRwavelet,OGFB,OGLM}, reconstruction \cite{GFTSSS}, classification \cite{GFTsampling}, clustering \cite{JFT,JFRFT,JLCT,Gclustering}, segmentation \cite{Segmentation}, estimation \cite{estimation1,estimation2}, detection \cite{detection}, semi-supervised learning \cite{Glearning}, graph learning \cite{graphlearning1,graphlearning2}, and the development of graph neural networks \cite{Gdeep,Gmachine,GCN}.

However, many graph-structured data also evolve over time. For instance, video sequences can be viewed as temporally evolving images, and user behaviors in social networks often vary dynamically. This motivates the need for a unified framework for analyzing joint time-vertex signals. Pioneering works \cite{JFT,HGFT} integrated discrete-time processing with GSP, laying the foundation for this direction. The joint time-vertex Fourier transform (JFT), which combines the DFT in the temporal domain and the GFT in the vertex domain, serves as a key tool for joint spectral analysis. Building on this framework, a variety of applications have emerged, including reconstruction of time-varying graph signals \cite{HGFRFT,Jreconstruction,Jsamp2,Jsamp3}, prediction of joint spectral-temporal information \cite{prediction}, and modeling of stationary dynamics on graphs \cite{evolution}. Spatio-temporal image and video restoration can be regarded as a special case of joint time-vertex signals, further highlighting the importance of developing joint transforms and filtering methods \cite{JFT,spatio-temporal}.

\IEEEpubidadjcol

To further enhance the expressive power of joint time-vertex signal processing, \textit{multiresolution analysis} plays an important role. Multiresolution analysis has long been effective in traditional analysis, compression, and denoising \cite{Multiresolution}, and its extension to the joint domain provides a powerful tool for decomposition and representation. In GSP, two-channel critically sampled wavelet filter banks have been proposed to realize multiresolution transforms \cite{biorth,PRwavelet}. For bipartite graphs, such filter banks avoid spectral folding during downsampling and ensure perfect reconstruction \cite{OGFB}. However, for arbitrary non-bipartite graphs, the original graph must be partitioned into bipartite subgraphs, each retaining all vertices but only a subset of edges, whose union reconstructs the original graph. Applying filter banks to each subgraph leads to multi-dimensional transforms, yet many original edges are not fully utilized in a single stage. To address this limitation, oversampled graph Laplacian matrices were introduced in \cite{OGFB,OGLM,OGFT}, enabling flexible oversampling through the addition of nodes and edges, and thereby supporting redundant signal representations. Notably, the filter bank concept has proven effective for graph-based decomposition and reconstruction and has inspired graph-based image denoising.

Image denoising aims to remove noise while preserving structural details, and GSP has attracted significant attention for this purpose \cite{GSIP}. Early studies employed graph Laplacian regularization to enhance smoothness while retaining edges, and subsequent extensions with manifold-based gradient regularizers and plug-and-play unrolled optimization approaches further improved denoising performance and flexibility \cite{GLRforID,GGLR,GLRUnrolling}. Multiresolution and reversible graph Fourier transforms were developed for denoising, as well as image compression and watermarking \cite{GFTimage1,GFTimage2}. In parallel, graph-based wavelet and filter bank methods have been proposed for image modeling. For example, graph wavelet filter banks provide sparse, edge-aware image representations, achieving superior performance over conventional wavelets in nonlinear approximation and denoising \cite{biorth,PRwavelet}. The oversampling concept further enhances these filter banks, improving robustness in image denoising and reconstruction \cite{OGFB,OGLM,OGFT}. These works demonstrate the strong potential of GSP in modern image and video denoising and restoration, motivating extensions to joint time-vertex graph methodologies.

Against this background, this paper develops a unified framework for joint time-vertex signal processing. The framework is built on an oversampled graph Laplacian matrix extended from the vertex domain to the temporal domain, supporting bipartite expansion and redundant multiresolution representations. Based on this operator, a two-channel joint time-vertex oversampled filter bank is designed to satisfy perfect reconstruction on bipartite graphs and realize oversampling in both time and vertex domains via a $K$-coloring strategy. The framework is applied to signal and image/video denoising, with experiments demonstrating improved decomposition and reconstruction performance compared to existing methods, highlighting its potential at the intersection of GSP and image processing.

To summarize, our main contributions include:
\begin{itemize}
	\item Introduce a joint time-vertex oversampled graph Laplacian matrixthat extends vertex domain oversampling to the temporal domain, preserving the original joint graph structure while enabling bipartite expansion and redundant multiresolution representations.
	\item Design a two-channel joint time-vertex oversampled filter bank that satisfies perfect reconstruction conditions on bipartite graphs and implements oversampling in both time and vertex domains using a $K$-coloring strategy.
	\item Apply the proposed framework to joint time-vertex signal and image/video denoising, demonstrating its effectiveness in decomposition and reconstruction, and outperforming existing methods in numerical experiments.
\end{itemize}

The remainder of the paper is organized as follows. Section II reviews joint time-vertex signal processing, two-channel oversampled graph filter banks, and oversampled graph Laplacians. Section III introduces the proposed time oversampled Laplacian and joint oversampled Laplacian, together with their properties. Section IV presents the design of critically sampled and oversampled joint filter banks. Section V discusses bipartite extensions in both time and vertex domains. Section VI reports experimental results on signal decomposition, reconstruction, and denoising, including image and video applications. Finally, Section VII concludes the paper.

		
		
\section{Preliminaries}

A brief overview of existing theories on joint time-vertex signal processing \cite{JFT}, graph filter banks, and oversampled graph Laplacian matrix \cite{OGFB,OGLM} is provided below.

\subsection{Joint Time-Vertex Signal Processing}
Let $\mathcal{G} = (\mathcal{V}, \mathcal{E}, \mathbf{A}_{\mathcal{G}})$ be a finite, undirected graph with no self-loops or multiple edges, where $\mathcal{V}$ and $\mathcal{E}$ denote the node and edge sets, respectively, and $|\mathcal{V}| = N$. The adjacency matrix $\mathbf{A}_{\mathcal{G}}$ has entries $a_{mn} = w_{mn}$ if nodes $m$ and $n$ are connected, and $0$ otherwise. The degree matrix $\mathbf{D}_{\mathcal{G}}$ is diagonal with $d_{mm} = \sum_n a_{mn}$. The graph Laplacian matrix is defined as $\mathbf{L}_{\mathcal{G}} = \mathbf{D}_{\mathcal{G}}-\mathbf{A}_{\mathcal{G}}$, and the symmetric normalized graph Laplacian is
\begin{equation}
	\bm{\mathcal{L}}_{\mathcal{G}} = \mathbf{D}^{-1/2}_{\mathcal{G}}\mathbf{L}_{\mathcal{G}}\mathbf{D}^{-1/2}_{\mathcal{G}},
\end{equation}
whose eigenvalues $0 = \lambda_{0} < \cdots \leq \lambda_{N-1} \leq 2$ define the diagonal matrix $\mathbf{\Lambda}_{\mathcal{G}}$. Let $\bm{u}_{\lambda_{\ell}}$ be the eigenvector of $\bm{\mathcal{L}}_{\mathcal{G}}$ for $\lambda_{\ell}$, and $\mathbf{U}_{\mathcal{G}} = [\bm{u}_{\lambda_0}, \dots, \bm{u}_{\lambda_{N-1}}]$.

Assume that a graph signal is sampled at $T$ regularly spaced discrete time instants. Let $\bm{x}_t \in \mathbb{R}^N$ denote the graph signal at time $t$, and define the time-varying matrix 
\[
\mathbf{X} = [\bm{x}_1, \bm{x}_2, \ldots, \bm{x}_T] \in \mathbb{R}^{N \times T}. 
\]
Both $\mathbf{X}$ and its vectorized form $\bm{x} = \text{vec}(\mathbf{X}) \in \mathbb{R}^{NT}$ are referred to as joint time-vertex signals \cite{JFT}. In $\mathbf{X}$, each column corresponds to the state of the graph at a specific time, while each row describes the temporal evolution of a single vertex. The GFT \cite{GFTlaplace,GFTadjacency1,GFTadjacency2} along the vertex dimension is defined as
$$
\text{GFT}\left\{ \mathbf{X} \right\} = \mathbf{U}^{\mathrm{H}}_{\mathcal{G}} \mathbf{X},
$$
and the DFT \cite{DFT} along the time dimension is given by
$$
\text{DFT}\left\{ \mathbf{X} \right\} = \mathbf{X} \mathbf{U}_{\mathcal{T}}^{\ast},
$$
where $\mathbf{U}_{\mathcal{T}}^{\ast}$ is the complex conjugate of the matrix $\mathbf{U}_{\mathcal{T}}(k, t) = \frac{1}{\sqrt{T}} \exp( \mathrm{i} \frac{2\pi(k - 1)t}{T}) $, for $k,t = 1,2,...,T$. Consequently, the JFT of $\mathbf{X}$ is expressed as
\[
\text{JFT}\left\{ \mathbf{X} \right\} = \mathbf{U}^{\mathrm{H}}_{\mathcal{G}} \mathbf{X} \mathbf{U}_{\mathcal{T}}^{\ast}.
\]
In vectorized form, the JFT becomes
\begin{equation}
	\bm{\hat{x}} = \text{JFT}\left\{ \bm{x} \right\} = \mathbf{U}^{\mathrm{H}}_{\mathcal{J}} \bm{x} = (\mathbf{U}^{\mathrm{H}}_{\mathcal{T}} \otimes \mathbf{U}^{\mathrm{H}}_{\mathcal{G}}) \bm{x}, \label{JFT}
\end{equation}
where $\otimes$ denotes the Kronecker product, and $\mathbf{U}_{\mathcal{J}}  = \mathbf{U}_{\mathcal{T}} \otimes \mathbf{U}_{\mathcal{G}}$ is the joint Fourier basis. The inverse JFT (IJFT) is given by $\bm{x} = \mathbf{U}_{\mathcal{J}}\bm{\hat{x}}$, which reconstructs the extended signal from its spectral representation.

Temporal smoothness is modeled using the normalized time Laplacian $\bm{\mathcal{L}}_{\mathcal{T}}$, defined on a ring graph with periodic boundary conditions. As each node in the temporal graph has degree 2, the normalized Laplacian is related to its unnormalized counterpart $\mathbf{L}_{\mathcal{T}}$ via $\bm{\mathcal{L}}_{\mathcal{T}} = \frac{1}{2} \mathbf{L}_{\mathcal{T}}$. The eigendecomposition of $\bm{\mathcal{L}}_{\mathcal{T}}$ is given by $\bm{\mathcal{L}}_{\mathcal{T}} = \mathbf{U}_{\mathcal{T}} \bm{\Lambda}_{\mathcal{T}} \mathbf{U}_{\mathcal{T}}^{\mathrm{H}}$, where $\bm{\Lambda}_{\mathcal{T}}(k,k) =\omega_k$ \cite{JFT}.

Using the mixed-product property of the Kronecker product, the joint time-vertex normalized Laplacian operator is defined as $\bm{\mathcal{L}}_{\mathcal{J}} = \bm{\mathcal{L}}_{\mathcal{T}} \otimes \mathbf{I}_N + \mathbf{I}_T \otimes \bm{\mathcal{L}}_{\mathcal{G}}$, which is jointly diagonalizable via
\begin{equation}
	\begin{aligned}\bm{\mathcal{L}}_{\mathcal{J}} =&(\mathbf{U}_{\mathcal{T}} \otimes \mathbf{U}_{\mathcal{G}}) (\bm{\Lambda}_{\mathcal{T}} \times \bm{\Lambda}_{\mathcal{G}}) (\mathbf{U}_{\mathcal{T}} \otimes \mathbf{U}_{\mathcal{G}})^{\mathrm{H}} \\ =&\mathbf{U}_{\mathcal{J}} \bm{\Lambda}_{\mathcal{J}} \mathbf{U}^{\mathrm{H}}_{\mathcal{J}},\end{aligned}  \label{LJ}
\end{equation}
where $\times$ denotes the sum of eigenvalues over all $(\omega_k, \lambda_{\ell})$ pairs across the time and vertex domains.

\subsection{Two-Channel Oversampled Graph Filter Banks}
A two-channel graph filter bank extends classical signal decomposition to graph-structured data. Consider a bipartite graph $\mathcal{G} = \{\mathcal{L}_\mathcal{G}, \mathcal{H}_\mathcal{G}, \mathcal{E}\}$ with $N$ nodes partitioned into low-pass and high-pass vertex sets \cite{OGFB}. The downsampling-then-upsampling operator for each channel is defined as
\[
\mathbf{D}^{du}_{\mathcal{L}_\mathcal{G}} = \frac{1}{2}(\mathbf{I} + \mathbf{C}_\mathcal{G}), \quad \mathbf{D}^{du}_{\mathcal{H}_\mathcal{G}} = \frac{1}{2}(\mathbf{I} - \mathbf{C}_\mathcal{G}),
\]
where $\mathbf{C}_\mathcal{G}$ is a diagonal matrix with entries
\begin{equation}
	(\mathbf{C}_\mathcal{G})_{mm} = 
	\begin{cases}
		+1, & \text{if } m \in \mathcal{L}_\mathcal{G}, \\
		-1, & \text{if } m \in \mathcal{H}_\mathcal{G}.
	\end{cases} \label{C}
\end{equation}

\begin{figure}[t]
	\begin{center}
		\begin{minipage}[t]{0.9\linewidth}
			\centering
			\includegraphics[width=\linewidth]{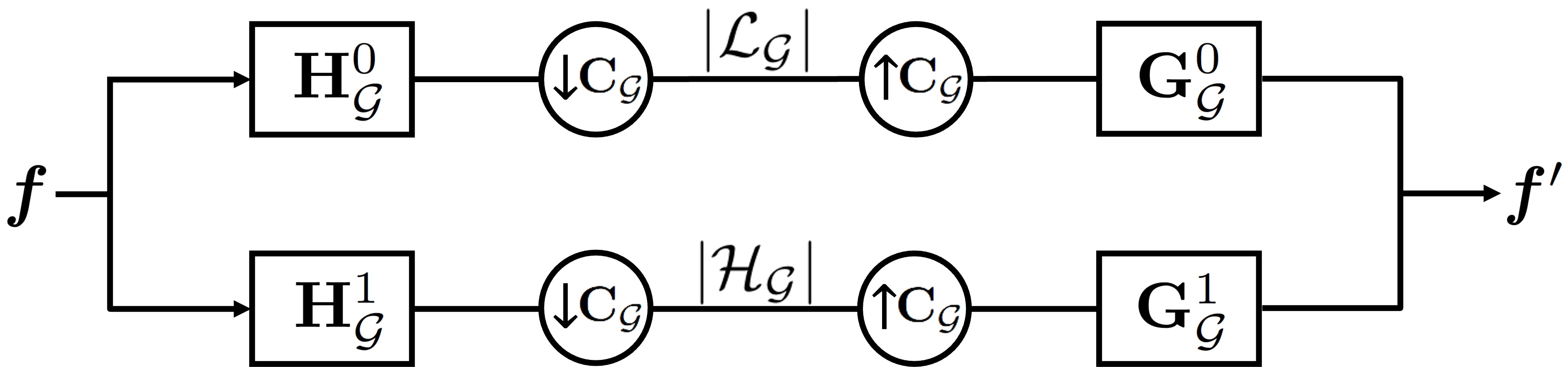}
		\end{minipage}
	\end{center}
	\caption{Two-channel filter bank based on critical graph sampling.}
	\vspace*{-3pt}
	\label{fig01}
\end{figure}

After applying the \textit{analysis filters} $\mathbf{H}^0_{\mathcal{G}}$ and $\mathbf{H}^1_{\mathcal{G}}$ to the input graph signal $\bm{f} \in \mathbb{R}^N$, the partial reconstructions from each channel are
\begin{equation}
	\bm{f}'_k =
	\begin{cases}
		\frac{1}{2} \mathbf{G}^k_{\mathcal{G}} (\mathbf{I} - \mathbf{C}) \mathbf{H}^k_{\mathcal{G}} \bm{f}, & k = 0, \\
		\frac{1}{2} \mathbf{G}^k_{\mathcal{G}} (\mathbf{I} + \mathbf{C}) \mathbf{H}^k_{\mathcal{G}} \bm{f}, & k = 1,
	\end{cases} \label{HG}
\end{equation}
where $\mathbf{H}^k_{\mathcal{G}} = \mathbf{U}_{\mathcal{G}} h^k(\bm{\Lambda}_{\mathcal{G}}) \mathbf{U}_{\mathcal{G}}^{\mathrm{H}}$ and $\mathbf{G}^k_{\mathcal{G}} = \mathbf{U}_{\mathcal{G}} g^k(\bm{\Lambda}_{\mathcal{G}}) \mathbf{U}_{\mathcal{G}}^{\mathrm{H}}$ denote the spectral-domain analysis and \textit{synthesis filters}, defined using kernel functions $h^k(\cdot)$ and $g^k(\cdot)$. This formulation aligns with the localized (bandlimited) graph operators defined in \cite{GUncertainty}. The reconstruction framework is shown in Fig.~\ref{fig01}.

Perfect reconstruction is achieved if the following spectral conditions are satisfied
\begin{equation}
	\begin{aligned}g^0(\lambda )h^0(\lambda )+g^1(\lambda )h^1(\lambda )=&2,\\ g^0(\lambda )h^0(2-\lambda )-g^1(\lambda )h^1(2-\lambda )=&0.\end{aligned}  \label{hg}
\end{equation}
These can be satisfied by a symmetric filter design, $g^0(\lambda) = h^1(2 - \lambda)$, and $g^1(\lambda) = h^0(2 - \lambda)$.

\subsection{Oversampled Graph Laplacian Matrix}
\label{secII.C}

To enhance flexibility in filter bank design, an oversampled bipartite graph $\widetilde{\mathcal{G}} = \{\widetilde{\mathcal{L}}_\mathcal{G}, \widetilde{\mathcal{H}}_\mathcal{G}, \widetilde{\mathcal{E}}\}$ is constructed such that the original bipartite sets satisfy $\mathcal{L}_\mathcal{G} \subseteq \widetilde{\mathcal{L}}_\mathcal{G}$ and $\mathcal{H}_\mathcal{G} \subseteq \widetilde{\mathcal{H}}_\mathcal{G}$. The number of nodes is expanded from $N_0$ to $N_1$, with $|\widetilde{\mathcal{L}}_\mathcal{G}| + |\widetilde{\mathcal{H}}_\mathcal{G}| = N_1 > N_0$ \cite{OGLM}. Let $\mathbf{A}_{\mathcal{G}_0}$ denote the adjacency matrix of the original graph $\mathcal{G}$. The adjacency matrix of the oversampled graph $\widetilde{\mathbf{A}}_{\mathcal{G}}$ is defined as
\begin{equation}
	\widetilde{\mathbf{A}}_{\mathcal{G}} = \begin{bmatrix}
		\mathbf{A}_{\mathcal{G}_0} & \mathbf{A}_{\mathcal{G}_{01}} \\
		\mathbf{A}_{\mathcal{G}_{01}}^\top & \mathbf{0}_{N_1 -N_0}
	\end{bmatrix},
\end{equation}
where $\mathbf{A}_{\mathcal{G}_{01}}$ encodes connections between original and added nodes. The degree matrix $\widetilde{\mathbf{D}}_{\mathcal{G}}$ is derived from $\widetilde{\mathbf{A}}_{\mathcal{G}}$, and the Laplacian is given by $\widetilde{\mathbf{L}}_{\mathcal{G}} = \widetilde{\mathbf{D}}_{\mathcal{G}} - \widetilde{\mathbf{A}}_{\mathcal{G}}$. The symmetric normalized Laplacian of the oversampled graph is then defined as
\begin{equation}
	\widetilde{\bm{\mathcal{L}}}_{\mathcal{G}} = \widetilde{\mathbf{D}}^{-1/2}_{\mathcal{G}} \widetilde{\mathbf{L}}_{\mathcal{G}} \widetilde{\mathbf{D}}^{-1/2}_{\mathcal{G}}.
\end{equation}

\section{Joint Time-Vertex Oversampled Graph Laplacian Matrix}
\label{sec: JOGLM}
In this section, we characterize the structure of the joint time-vertex oversampled graph Laplacian and the associated oversampled graph signals, based on the oversampled graph Laplacian matrices \cite{OGLM}. We also present key properties and techniques of the joint Laplacian. Prior to applying the joint graph filter banks, the original graph $\mathcal{G}$ is extended to its oversampled version $\widetilde{\mathcal{G}} = \{ \widetilde{\mathcal{L}}_\mathcal{G}, \widetilde{\mathcal{H}}_\mathcal{G}\}$; similarly, the original temporal bipartite graph $\mathcal{T} = \{ \mathcal{L}_\mathcal{T}, \mathcal{H}_\mathcal{T} \}$ is extended to $\widetilde{\mathcal{T}} = \{ \widetilde{\mathcal{L}}_\mathcal{T}, \widetilde{\mathcal{H}}_\mathcal{T}\}$. The joint spectral filtering is then performed using the constructed joint time-vertex oversampled Laplacian.

\subsection{Oversampled Time Laplacian Matrix}
The vertex domain oversampling is introduced by constructing an augmented bipartite graph $\widetilde{\mathcal{G}}$, which extends the original graph $\mathcal{G}_0$ with additional nodes to enable more sampling channels. Analogously, we now consider oversampling in the time domain by expanding the temporal graph structure.

Let the time domain be modeled as a ring graph $\mathcal{T}=\{\mathcal{V}_{\mathcal{T}}, \mathcal{E}_{\mathcal{T}}\}$, where $\mathcal{V}_{\mathcal{T}}$ denotes the set of time nodes, and $\mathcal{E}_{\mathcal{T}}$ denotes the set of edges capturing pairwise temporal dependencies. This cyclic structure encodes periodic relationships between time instances. We define its oversampled version as $\widetilde{\mathcal{T}} = \{\widetilde{\mathcal{L}}_\mathcal{T}, \widetilde{\mathcal{H}}_\mathcal{T}, \widetilde{\mathcal{E}}_{\mathcal{T}} \}$, where $\widetilde{\mathcal{L}}_\mathcal{T}$ and $\widetilde{\mathcal{H}}_\mathcal{T}$ denote the low-pass and high-pass node sets in the augmented temporal graph, respectively. The total number of time nodes increases from $T_0$ to $T_1$, i.e., $|\widetilde{\mathcal{L}}_\mathcal{T}| + |\widetilde{\mathcal{H}}_\mathcal{T}| = T_1 > T_0$, mirroring the expansion strategy employed in the vertex domain \cite{OGLM}.

This temporal oversampling can be interpreted as inserting auxiliary nodes into the periodic ring graph, thereby increasing the temporal resolution and allowing for additional sampling branches in the time domain.

Based on the adjacency matrix $\mathbf{A}_{\mathcal{T}}$ of the original temporal ring graph, we define the adjacency matrix of the oversampled temporal graph as
\begin{equation}
	\widetilde{\mathbf{A}}_{\mathcal{T}} =
	\begin{bmatrix}
		\mathbf{A}_{\mathcal{T}_0} & \mathbf{A}_{\mathcal{T}_{01}} \\
		\mathbf{A}_{\mathcal{T}_{01}}^\top & \mathbf{0}_{T_1 - T_0}
	\end{bmatrix},
\end{equation}
where $\mathbf{A}_{\mathcal{T}_0} \in \mathbb{R}^{T_0 \times T_0}$ is the original temporal adjacency matrix, and $\mathbf{A}_{\mathcal{T}_{01}} \in \mathbb{R}^{T_0 \times (T_1 - T_0)}$ encodes the connections between the original and newly added time nodes.

The degree matrix $\widetilde{\mathbf{D}}_{\mathcal{T}}$ is computed from $\widetilde{\mathbf{A}}_{\mathcal{T}}$, and the symmetric normalized Laplacian of the oversampled temporal graph is then defined as
\begin{equation}
\widetilde{\bm{\mathcal{L}}}_{\mathcal{T}} = \widetilde{\mathbf{D}}_{\mathcal{T}}^{-1/2} (\widetilde{\mathbf{D}}_{\mathcal{T}} - \widetilde{\mathbf{A}}_{\mathcal{T}}) \widetilde{\mathbf{D}}_{\mathcal{T}}^{-1/2}.
\end{equation}

This operator retains the desirable spectral properties of the original time Laplacian and simultaneously plays a central role in constructing the joint time-vertex oversampled Laplacian.

\subsection{Oversampled Time-Vertex Laplacian Matrix}
Building upon the oversampled graph constructions in both the time and vertex domains, we now propose a joint time-vertex oversampled graph Laplacian to support unified modeling and filter bank design across both dimensions.

Let the original time–vertex signal $\mathbf{X}$ be defined on the product graph $\mathcal{T} \otimes \mathcal{G}$ \cite{JFT}. By inserting auxiliary nodes in both time and vertex domains, we construct an extended joint graph $\widetilde{\mathcal{T}} \otimes \widetilde{\mathcal{G}}$ with enlarged sizes $T_1 > T_0$ and $N_1 > N_0$ in each respective domain. Under this structure, the oversampled graph $\widetilde{\mathcal{G}}$ is replicated $T_1$ times, corresponding to each oversampled time step, while inter-layer temporal connections are established based on the edge structure of $\widetilde{\mathcal{T}}$ to connect the same node across adjacent time layers.

\textit{Definition 1:} Let $\widetilde{\bm{\mathcal{L}}}_{\mathcal{T}}$ and $\widetilde{\bm{\mathcal{L}}}_{\mathcal{G}}$ denote the symmetric normalized Laplacian matrices of the oversampled temporal and vertex graphs, respectively. Then, the joint time-vertex Laplacian over the oversampled graph is defined via the Kronecker sum as
\begin{equation}
	\widetilde{\bm{\mathcal{L}}}_{\mathcal{J}} = \widetilde{\bm{\mathcal{L}}}_{\mathcal{T}} \times \widetilde{\bm{\mathcal{L}}}_{\mathcal{G}} = \widetilde{\bm{\mathcal{L}}}_{\mathcal{T}} \otimes \mathbf{I}_{N_1} + \mathbf{I}_{T_1} \otimes \widetilde{\bm{\mathcal{L}}}_{\mathcal{G}}, \label{widetildeLJ}
\end{equation}
where $\widetilde{\bm{\mathcal{L}}}_{\mathcal{J}}$ is of dimension $N_1 T_1 \times N_1 T_1$.

This construction guarantees that the spectral structure of the joint Laplacian remains separable. Although Eq.~\eqref{widetildeLJ} is structurally similar to the separable joint Laplacian operators \cite{JFT} defined by Eq.~\eqref{LJ}, our framework does not operate directly on the original product graph. Instead, both the time and vertex domains are first extended through structured oversampling. Therefore, the resulting joint Laplacian is defined on an augmented graph that supports a redundant two-channel filter bank design, rather than a purely separable joint filtering framework.

Let $\widetilde{\mathbf{U}}_{\mathcal{T}}$ and $\widetilde{\mathbf{U}}_{\mathcal{G}}$ denote the eigenvector matrices of the oversampled temporal and vertex Laplacians, with corresponding eigenvalue diagonal matrices $\widetilde{\bm{\Lambda}}_{\mathcal{T}}$ and $\widetilde{\bm{\Lambda}}_{\mathcal{G}}$. The spectral decomposition of $\widetilde{\bm{\mathcal{L}}}_{\mathcal{J}}$ is then given by
\begin{equation}
\begin{aligned}
	\widetilde{\bm{\mathcal{L}}}_{\mathcal{J}} =
	& (\widetilde{\mathbf{U}}_{\mathcal{T}} \widetilde{\bm{\Lambda}}_{\mathcal{T}} \widetilde{\mathbf{U}}_{\mathcal{T}}^\mathrm{H}) \otimes \mathbf{I}_{N_1} + \mathbf{I}_{T_1} \otimes (\widetilde{\mathbf{U}}_{\mathcal{G}} \widetilde{\bm{\Lambda}}_{\mathcal{G}} \widetilde{\mathbf{U}}_{\mathcal{G}}^\mathrm{H}) \\
	=& (\widetilde{\mathbf{U}}_{\mathcal{T}} \otimes \widetilde{\mathbf{U}}_{\mathcal{G}}) \left( \widetilde{\bm{\Lambda}}_{\mathcal{T}} \otimes \mathbf{I}_{N_1} + \mathbf{I}_{T_1} \otimes \widetilde{\bm{\Lambda}}_{\mathcal{G}} \right) (\widetilde{\mathbf{U}}_{\mathcal{T}} \otimes \widetilde{\mathbf{U}}_{\mathcal{G}})^\mathrm{H}\\
	=& \widetilde{\mathbf{U}}_{\mathcal{J}} \widetilde{\bm{\Lambda}}_{\mathcal{J}} \widetilde{\mathbf{U}}^{\mathrm{H}}_{\mathcal{J}}.
\end{aligned}
\end{equation}
This formulation reveals that the joint spectrum is composed by summing all possible pairs of temporal and vertex eigenvalues, while the joint eigenvectors are constructed as the Kronecker product of their respective bases. This separable spectral structure facilitates efficient implementation of the proposed joint oversampled filter bank, while the redundancy introduced by oversampling fundamentally distinguishes our framework from conventional separable joint filtering methods.

To align with the joint time-vertex oversampled graph structure $\widetilde{\mathcal{T}} \otimes \widetilde{\mathcal{G}}$, the original time-vertex signal must be accordingly extended. Let the original signal be denoted as $\mathbf{X}_{0,0} \in \mathbb{R}^{N_0 \times T_0}$, defined on the original product graph $\mathcal{T} \otimes \mathcal{G}$, where $T_0$ and $N_0$ are the numbers of temporal and vertex nodes, respectively.

After inserting auxiliary nodes in both the temporal and vertex domains, we obtain the extended joint graph $\widetilde{\mathcal{T}} \otimes \widetilde{\mathcal{G}}$. Correspondingly, the joint oversampled graph signal is denoted by $\widetilde{\mathbf{X}} \in \mathbb{R}^{N_1 \times T_1}$, and can be organized in a block-matrix form as follows
\begin{equation}
\widetilde{\mathbf{X}} = 
\begin{bmatrix}
	\mathbf{X}_{0,0} & \mathbf{X}_{0,1} \\
	\mathbf{X}_{1,0} & \mathbf{X}_{1,1}
\end{bmatrix}, \label{widetildeX}
\end{equation}
where $\mathbf{X}_{0,1} \in \mathbb{R}^{(N_1 - N_0) \times T_0}$ corresponds to the newly introduced vertex nodes at the original time steps, $\mathbf{X}_{1,0} \in \mathbb{R}^{N_0 \times (T_1 - T_0)}$ represents the original vertex nodes at the newly inserted time steps, and $\mathbf{X}_{1,1} \in \mathbb{R}^{(N_1 - N_0) \times (T_1 - T_0)}$ denotes the signal defined over both the new temporal and new vertex nodes.

\textit{Remark 1:} The values of these additional components may be initialized as zeros to facilitate theoretical analysis, estimated through interpolation methods, or treated as learnable parameters within optimization or learning frameworks, depending on the application context.

To quantify the variation of signals over the joint oversampled graph, we introduce a set of gradient-induced norms defined over three domains: time, vertex, and joint \cite{JFT}. 

\textit{Definition 2:}  Let $\widetilde{\mathbf{X}}$ be the joint oversampled signal and $\widetilde{\bm{x}} = \text{vec}(\widetilde{\mathbf{X}})$ its vectorized form. The gradient operators are defined as follows
\begin{itemize}
\item Time gradient: $\nabla_{\mathcal{T}} \widetilde{\mathbf{X}} := \widetilde{\mathbf{X}} \widetilde{\bm{\mathcal{L}}}_{\mathcal{T}}^{1/2}$; 
\item Vertex gradient: $\nabla_{\mathcal{G}} \widetilde{\mathbf{X}} := \widetilde{\bm{\mathcal{L}}}_{\mathcal{G}}^{1/2} \widetilde{\mathbf{X}}$; 
\item Joint time-vertex gradient: $\nabla_{\mathcal{J}} \widetilde{\bm{x}} := \widetilde{\bm{\mathcal{L}}}_{\mathcal{J}}^{1/2} \widetilde{\bm{x}}$.
\end{itemize}

We begin with the $\ell_2$-norm of the joint time-vertex gradient, which characterizes the overall smoothness or quadratic variation of the signal,
\[
\begin{aligned}
\| \nabla_{\mathcal{J}} \widetilde{\bm{x}} \|_2^2 &= \widetilde{\bm{x}}^\top \widetilde{\bm{\mathcal{L}}}_{\mathcal{J}} \widetilde{\bm{x}} = \| \nabla_{\mathcal{G}} \widetilde{\mathbf{X}} \|_F^2 + \| \nabla_{\mathcal{T}} \widetilde{\mathbf{X}} \|_F^2 \\
& =\operatorname{tr}(\widetilde{\mathbf{X}}^\top \widetilde{\bm{\mathcal{L}}}_{\mathcal{G}} \widetilde{\mathbf{X}}) + \operatorname{tr}(\widetilde{\mathbf{X}} \widetilde{\bm{\mathcal{L}}}_{\mathcal{T}} \widetilde{\mathbf{X}}^\top),
\end{aligned}
\]
where $\| \cdot \|_F$ is the Frobenius norm and $\operatorname{tr}(\cdot)$ is the trace operator, reflecting contributions from both vertex and time domains.
 
To better handle non-smooth structures or piecewise-constant signals, we further consider the $\ell_1$-norm of the joint gradient, also referred to as the total variation norm
\[
\| \nabla_{\mathcal{J}} \widetilde{\bm{x}} \|_1 = \| \text{vec}(\nabla_{\mathcal{G}} \widetilde{\mathbf{X}}) \|_1 + \| \text{vec}(\nabla_{\mathcal{T}} \widetilde{\mathbf{X}}) \|_1.
\]

Moreover, a generalized mixed-norm formulation can be adopted to flexibly balance sensitivity across the two domains
\[
\mathcal{N}_{p,q}(\widetilde{\mathbf{X}}) = \mu_{\mathcal{G}} \| \text{vec}(\nabla_{\mathcal{G}} \widetilde{\mathbf{X}}) \|_p^p + \mu_{\mathcal{T}} \| \text{vec}(\nabla_{\mathcal{T}} \widetilde{\mathbf{X}}) \|_q^q,
\]
where $\mu_{\mathcal{G}}, \mu_{\mathcal{T}} \geq 0$ are domain-specific weights and $p, q \geq 1$ are norm parameters. This formulation is particularly well-suited to modeling heterogeneous smoothness patterns, such as signals that are highly irregular on the graph structure but relatively smooth over time.

Building upon the eigendecomposition of the oversampled operators, the joint Fourier transform can be naturally extended to the oversampled setting. Specifically, for the extended graph signal $\widetilde{\mathbf{X}}$, the JFT is defined as
\begin{equation}
\mathcal{F}_{\widetilde{\mathbf{X}}}=\text{JFT}\left\{ \widetilde{\mathbf{X}} \right\} = \widetilde{\mathbf{U}}_{\mathcal{G}}^{\mathrm{H}} \widetilde{\mathbf{X}} \widetilde{\mathbf{U}}^{\ast}_{\mathcal{T}},
\end{equation}
in vectorized form, the transform becomes
$$
\mathcal{F}_{\widetilde{\bm{x}}}=\text{JFT}\left\{ \widetilde{\bm{x}} \right\} = \widetilde{\mathbf{U}}_{\mathcal{J}}^{\mathrm{H}} \widetilde{\bm{x}},
$$
where $\widetilde{\mathbf{U}}_{\mathcal{J}} = \widetilde{\mathbf{U}}_{\mathcal{T}} \otimes \widetilde{\mathbf{U}}_{\mathcal{G}}$, consistent with the definition in Eq.~\eqref{JFT}, and the transformed coefficients $\mathcal{F}_{\widetilde{\bm{x}}}$ lie in the joint spectral domain.

\textit{Remark 2:} Unlike conventional joint filtering on the original product graph, the proposed transform underpins a redundant joint time-vertex filter bank, enabling multiresolution representation and enhanced reconstruction stability.

\section{Design of Two-Channel Joint Time-Vertex Oversampled Graph Filter Banks}
\label{JOGFB}
In this section, we extend the critically sampled graph filter banks \cite{OGFB} from the vertex domain to the joint time-vertex domain, and further develop its natural generalization to the joint oversampled setting.

\subsection{Joint Critically Sampled Graph Filter Banks}
The time-vertex graph $\mathcal{J} = \mathcal{T} \otimes \mathcal{G}$ enables joint processing along both domains. The input signal $\mathbf{X}$ is first filtered in the vertex and temporal domains using analysis filters $\mathbf{H}^k_{\mathcal{G}}$ and $\mathbf{H}^k_{\mathcal{T}}$, where $\mathbf{H}^k_{\mathcal{G}}$ is defined as Eq. \eqref{HG}, and another is 
\begin{equation}
	\mathbf{H}^k_{\mathcal{T}} = \mathbf{U}_{\mathcal{T}} h^k(\bm{\Lambda}_{\mathcal{T}}) \mathbf{U}_{\mathcal{T}}^{\mathrm{H}},
\end{equation}
where $k = 0$ or $1$, thus the joint filtering result is given by
\begin{equation}
	\mathbf{X}^k_{\mathcal{J}} = \mathbf{H}^k_{\mathcal{G}} \mathbf{X} (\mathbf{H}^k_{\mathcal{T}})^*.
\end{equation}

The sampling diagonal matrices $\mathbf{C}_{\mathcal{G}}$ and $\mathbf{C}_{\mathcal{T}}$, consistent with the definition in Eq. \eqref{C} are employed to perform structured downsampling along the vertex and time axes, respectively.

The reconstruction process uses synthesis filters $\mathbf{G}^k_{\mathcal{G}}$ and $\mathbf{G}^k_{\mathcal{T}}$, which share the same eigenspaces as the analysis filters but are modulated by synthesis kernels $g^k(\cdot)$,
$$
\begin{cases}
	\mathbf{G}^k_{\mathcal{G}} = \mathbf{U}_{\mathcal{G}} g^k(\bm{\Lambda}_{\mathcal{G}}) \mathbf{U}_{\mathcal{G}}^{\mathrm{H}}, \\
	\mathbf{G}^k_{\mathcal{T}} = \mathbf{U}_{\mathcal{T}} g^k(\bm{\Lambda}_{\mathcal{T}}) \mathbf{U}_{\mathcal{T}}^{\mathrm{H}}.
\end{cases}
$$
Therefore, the full joint reconstruction is then given by
\begin{equation}
	\begin{aligned}
		\mathbf{X}'_{\mathcal{J}} = 
		&\frac{1}{2} \mathbf{G}^0_{\mathcal{G}} (\mathbf{I} - \mathbf{C}_{\mathcal{G}}) \mathbf{H}^0_{\mathcal{G}} \, \mathbf{X} \, (\mathbf{H}^0_{\mathcal{T}})^* (\mathbf{I} - \mathbf{C}_{\mathcal{T}}) (\mathbf{G}^0_{\mathcal{T}})^* \\
		&+\frac{1}{2} \mathbf{G}^1_{\mathcal{G}} (\mathbf{I} + \mathbf{C}_{\mathcal{G}}) \mathbf{H}^1_{\mathcal{G}} \, \mathbf{X} \, (\mathbf{H}^1_{\mathcal{T}})^* (\mathbf{I} + \mathbf{C}_{\mathcal{T}}) (\mathbf{G}^1_{\mathcal{T}})^*.
	\end{aligned}
\end{equation}

The filter bank achieves perfect reconstruction if the following spectral conditions are satisfied for all eigenvalues $\lambda$ and $\omega$,
\begin{equation}
	\begin{cases}
		g^0(\lambda) h^0(\lambda) + g^1(\lambda) h^1(\lambda) = 2, \\
		g^0(\lambda) h^0(2 - \lambda) - g^1(\lambda) h^1(2 - \lambda) = 0,\\
		g^0(\omega) h^0(\omega) + g^1(\omega) h^1(\omega) = 2, \\
		g^0(\omega) h^0(2 - \omega) - g^1(\omega) h^1(2 - \omega) = 0.
	\end{cases} \label{condition}
\end{equation}
These can be satisfied via symmetric filter designs as Eq. \eqref{hg}, and similarly in the temporal domain.

\subsection{Joint Oversampled Graph Filter Banks}
To enable flexible signal modeling and transformation in the joint time-vertex domain, we introduce two-channel oversampled graph filter banks. Unlike the critically sampled case \cite{PRwavelet}, the oversampled structure preserves all nodes in both time and vertex domains, allowing richer representation and enabling integration into learning frameworks \cite{OGLM}.

The normalized Laplacian of a bipartite graph has a spectrum symmetric about 1, so if $\lambda$ is an eigenvalue, $2-\lambda$ is also an eigenvalue (Eq.~\eqref{hg}). This enables quadrature mirror kernels that cancel aliasing \cite{OGFB}. In the joint time-vertex setting, separability ensures that perfect reconstruction in each domain guarantees perfect reconstruction of the overall transform, extending classical graph filter bank to the joint domain.

The oversampled signal $\widetilde{\mathbf{X}} \in \mathbb{R}^{N_1 \times T_1}$ is defined on the joint graph $\widetilde{\mathcal{T}} \otimes \widetilde{\mathcal{G}}$. For each channel $k = 0, 1$, the input signal is filtered jointly in the vertex and time domains
\begin{equation}
	\widetilde{\mathbf{X}}^k_{\mathcal{J}} = \widetilde{\mathbf{H}}^k_{\mathcal{G}} \, \widetilde{\mathbf{X}} \, (\widetilde{\mathbf{H}}^k_{\mathcal{T}})^*,
\end{equation}
where the analysis filters are defined via the spectral decomposition,
\begin{equation}
	\begin{cases}
		\widetilde{\mathbf{H}}^k_{\mathcal{G}} = \widetilde{\mathbf{U}}_{\mathcal{G}} h^k(\widetilde{\bm{\Lambda}}_{\mathcal{G}}) \widetilde{\mathbf{U}}_{\mathcal{G}}^{\mathrm{H}},\\
		\widetilde{\mathbf{H}}^k_{\mathcal{T}} = \widetilde{\mathbf{U}}_{\mathcal{T}} h^k(\widetilde{\bm{\Lambda}}_{\mathcal{T}})  \widetilde{\mathbf{U}}_{\mathcal{T}}^{\mathrm{H}}.
	\end{cases}
\end{equation}

\textit{Remark 3:} No node is discarded in the oversampled setting; instead, the full structure is preserved to maintain signal completeness.

The reconstructed signal $\widetilde{\mathbf{X}}'_{\mathcal{J}}$ is obtained by combining the outputs of both channels, each processed by a synthesis filter and domain-specific selection matrix as follows
\begin{equation}
	\begin{aligned}
		\widetilde{\mathbf{X}}'_{\mathcal{J}} = & \frac{1}{2} \widetilde{\mathbf{G}}^0_{\mathcal{G}} (\mathbf{I}_{N_1} - \widetilde{\mathbf{C}}_{\mathcal{G}}) \widetilde{\mathbf{H}}^0_{\mathcal{G}} \, \widetilde{\mathbf{X}} \, (\widetilde{\mathbf{H}}^0_{\mathcal{T}})^* (\mathbf{I}_{T_1} - \widetilde{\mathbf{C}}_{\mathcal{T}}) (\widetilde{\mathbf{G}}^0_{\mathcal{T}})^* \\
		 &+ \frac{1}{2} \widetilde{\mathbf{G}}^1_{\mathcal{G}} (\mathbf{I}_{N_1} + \widetilde{\mathbf{C}}_{\mathcal{G}}) \widetilde{\mathbf{H}}^1_{\mathcal{G}} \, \widetilde{\mathbf{X}} \, (\widetilde{\mathbf{H}}^1_{\mathcal{T}})^* (\mathbf{I}_{T_1} + \widetilde{\mathbf{C}}_{\mathcal{T}}) (\widetilde{\mathbf{G}}^1_{\mathcal{T}})^*,
	\end{aligned}\label{OX'J}
\end{equation}
where $\widetilde{\mathbf{C}}_{\mathcal{G}}$ and $\widetilde{\mathbf{C}}_{\mathcal{T}}$ denote the oversampled sampling matrices in the vertex and time domains, respectively. Moreover, in Eq. \eqref{OX'J}, the synthesis filters are constructed in the same spectral form as the analysis filters, with the kernel function $g^k(\cdot)$ replacing $h^k(\cdot)$,
\begin{equation}
	\begin{cases}
		\widetilde{\mathbf{G}}^k_{\mathcal{G}} = \widetilde{\mathbf{U}}_{\mathcal{G}} \, g^k(\widetilde{\bm{\Lambda}}_{\mathcal{G}}) \, \widetilde{\mathbf{U}}_{\mathcal{G}}^{\mathrm{H}}, \\
		\widetilde{\mathbf{G}}^k_{\mathcal{T}} = \widetilde{\mathbf{U}}_{\mathcal{T}} \, g^k(\widetilde{\bm{\Lambda}}_{\mathcal{T}}) \, \widetilde{\mathbf{U}}_{\mathcal{T}}^{\mathrm{H}}.
	\end{cases}
\end{equation}

The filter bank guarantees perfect reconstruction if the spectral kernels satisfy the conditions in Eq.~\eqref{condition} for all eigenvalues $\lambda \in \mathrm{spec}(\widetilde{\bm{\mathcal{L}}}_{\mathcal{G}})$ and $\omega \in \mathrm{spec}(\widetilde{\bm{\mathcal{L}}}_{\mathcal{T}})$.

\begin{figure*}[t]
	\begin{center}
		\begin{minipage}[t]{1\linewidth}
			\centering
			\includegraphics[width=\linewidth]{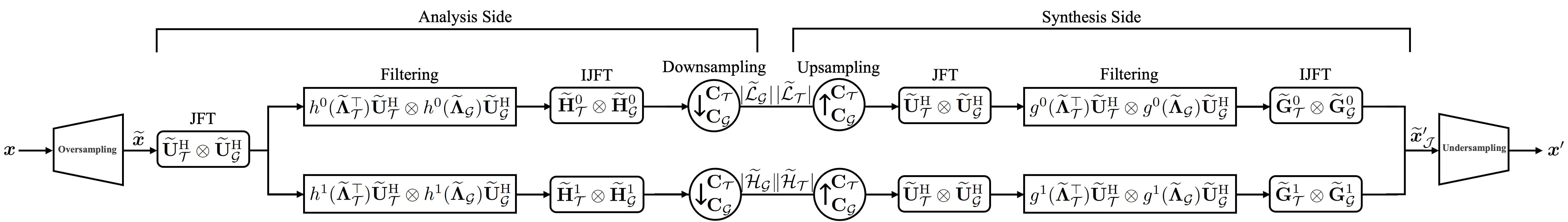}
		\end{minipage}
	\end{center}
	\caption{Framework of joint time-vertex oversampled graph transform.}
	\vspace*{-3pt}
	\label{fig02}
\end{figure*}

We further rewrite the joint reconstruction process in a vectorized form. By exploiting the properties of the Kronecker product, Eq.~\eqref{OX'J} can be equivalently expressed as
\[
\begin{aligned}
	\widetilde{\bm{x}}'_{\mathcal{J}} =\ & \frac{1}{2} \left( \widetilde{\mathbf{G}}^0_{\mathcal{T}} (\mathbf{I}_{T_1} - \widetilde{\mathbf{C}}_{\mathcal{T}}) \widetilde{\mathbf{H}}^0_{\mathcal{T}}  \otimes \widetilde{\mathbf{G}}^0_{\mathcal{G}} (\mathbf{I}_{N_1} - \widetilde{\mathbf{C}}_{\mathcal{G}}) \widetilde{\mathbf{H}}^0_{\mathcal{G}} \right) \widetilde{\bm{x}} \\
	+& \frac{1}{2} \left( \widetilde{\mathbf{G}}^1_{\mathcal{T}} (\mathbf{I}_{T_1} + \widetilde{\mathbf{C}}_{\mathcal{T}}) \widetilde{\mathbf{H}}^1_{\mathcal{T}} \otimes \widetilde{\mathbf{G}}^1_{\mathcal{G}} (\mathbf{I}_{N_1} + \widetilde{\mathbf{C}}_{\mathcal{G}}) \widetilde{\mathbf{H}}^1_{\mathcal{G}} \right) \widetilde{\bm{x}},
\end{aligned}
\]
where the analysis and synthesis filters $\widetilde{\mathbf{H}}^k_{\mathcal{T}}$, and $\widetilde{\mathbf{G}}^k_{\mathcal{T}}$ are Hermitian matrices, i.e., $(\cdot)^{\mathrm{H}} = (\cdot)$, which permits simplification of the conjugate transpose operations in the original expression. We denote the following composite operators for compactness
\begin{align*}
	\widetilde{\mathbf{Q}}^0_{\mathcal{T}} &= \widetilde{\mathbf{G}}^0_{\mathcal{T}} (\mathbf{I}_{T_1} - \widetilde{\mathbf{C}}_{\mathcal{T}}) \widetilde{\mathbf{H}}^0_{\mathcal{T}}, &
	\widetilde{\mathbf{Q}}^0_{\mathcal{G}} &= \widetilde{\mathbf{G}}^0_{\mathcal{G}} (\mathbf{I}_{N_1} - \widetilde{\mathbf{C}}_{\mathcal{G}}) \widetilde{\mathbf{H}}^0_{\mathcal{G}}, \\
	\widetilde{\mathbf{Q}}^1_{\mathcal{T}} &= \widetilde{\mathbf{G}}^1_{\mathcal{T}} (\mathbf{I}_{T_1} + \widetilde{\mathbf{C}}_{\mathcal{T}}) \widetilde{\mathbf{H}}^1_{\mathcal{T}}, &
	\widetilde{\mathbf{Q}}^1_{\mathcal{G}} &= \widetilde{\mathbf{G}}^1_{\mathcal{G}} (\mathbf{I}_{N_1} + \widetilde{\mathbf{C}}_{\mathcal{G}}) \widetilde{\mathbf{H}}^1_{\mathcal{G}}.
\end{align*}
The final reconstruction formula becomes
\begin{equation}
	\widetilde{\bm{x}}'_{\mathcal{J}} =
	\frac{1}{2} \left( \widetilde{\mathbf{Q}}^0_{\mathcal{T}} \otimes \widetilde{\mathbf{Q}}^0_{\mathcal{G}} \right) \widetilde{\bm{x}} +
	\frac{1}{2} \left( \widetilde{\mathbf{Q}}^1_{\mathcal{T}} \otimes \widetilde{\mathbf{Q}}^1_{\mathcal{G}} \right) \widetilde{\bm{x}}.
\end{equation}

These constraints ensure that both vertex and temporal dimensions are perfectly reconstructed under joint oversampled filtering. The overall architecture is illustrated in Fig.~\ref{fig02}.

\section{Efficient Joint Time-Vertex Oversampling Graph Extension Methods}
To mitigate the potential performance degradation caused by suboptimal oversampling designs, it is essential to impose a bipartite structure on general arbitrary graphs. In the vertex domain, graph coloring techniques such as Three-colorable or more generally $K$-colorable are among the most widely used approaches for constructing bipartite graphs \cite{biorth,PRwavelet,OGLM,OGFB}. Motivated by this, the concept of graph extension can be naturally extended to the joint time-vertex domain. By observing that the time-domain Laplacian coincides with the Laplacian of an undirected ring graph, we first propose an extension scheme in the time domain. This is then systematically combined with existing vertex-domain coloring methods to develop a unified graph extension framework in the joint domain.

\subsection{Construction of Time Bipartite Graph}
Within the joint time-vertex signal processing framework, we propose a method for constructing a bipartite graph structure in the time domain. The key idea is to model the temporal sequence as a ring graph and employ graph coloring theory to design an oversampled bipartite structure. The detailed construction is presented in Appendix \ref{appendix A}.

We model the time axis as a 1D ring graph of length $T$. When $T$ is even, i.e., $T = 2t$, the ring is a natural bipartite graph and can be directly partitioned into two disjoint color classes. For odd $T = 2t + 1$, however, the ring graph is not bipartite. According to graph coloring theory, such a ring is always three-colorable \cite{OGLM}. Let the three color classes be denoted by $\mathcal{S}_{\mathcal{T}_1}$, $\mathcal{S}_{\mathcal{T}_2}$, and $\mathcal{S}_{\mathcal{T}_3}$. Suppose that $\mathcal{S}_{\mathcal{T}_1}$ and $\mathcal{S}_{\mathcal{T}_3}$ each contain $t$ time points, while $\mathcal{S}_{\mathcal{T}_2}$ contains a single node.

The foundation bipartite graph is constructed using the partition $\mathcal{S}_{\mathcal{T}_1} \cup \mathcal{S}_{\mathcal{T}_2}$ and $\mathcal{S}_{\mathcal{T}_3}$, with edges retained only between these two sets. To restore the internal edges within $\mathcal{S}_{\mathcal{T}_1} \cup \mathcal{S}_{\mathcal{T}_2}$ while preserving bipartiteness, the nodes are duplicated in this set, generating $\mathcal{S}'_{\mathcal{T}_1}$ and $\mathcal{S}'_{\mathcal{T}_2}$, and assign them to the opposite bipartite layer.

\begin{figure}[h!]
	\centering
	\includegraphics[width=\linewidth]{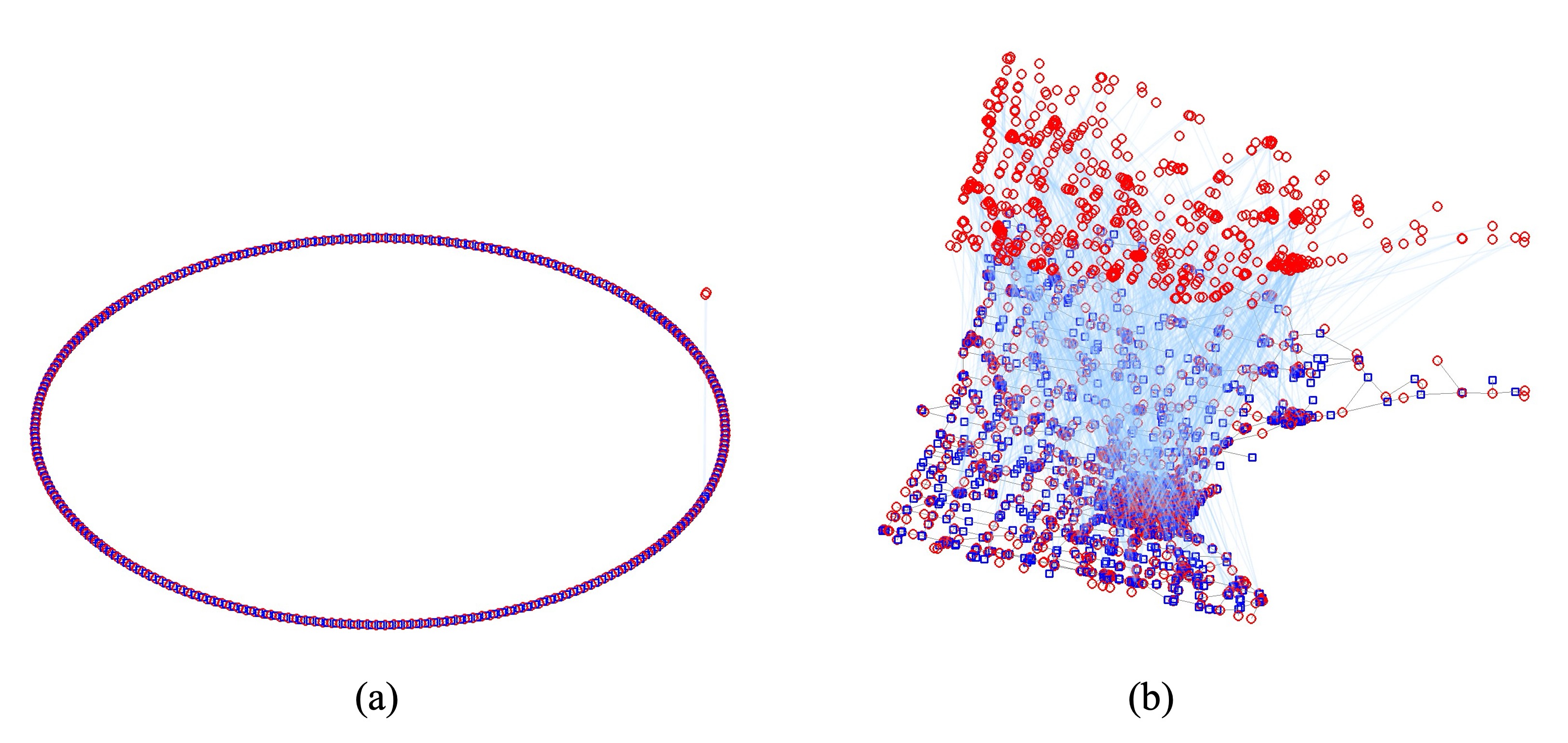}
	\vspace*{-15pt}
	\caption{Oversampled bipartite graph construction in the time and vertex domains. (a) The ring graph ($T=487$). (b) The Minnesota road network. The blue squares and red circles represent the low-pass and high-pass sets, respectively, and the light-blue edges indicate the newly added connections.}
	\label{fig03}
\end{figure}

The resulting oversampled bipartite graph has node sets $\widetilde{\mathcal{L}}_{\mathcal{T}} = \mathcal{S}_{\mathcal{T}_1} \cup \mathcal{S}_{\mathcal{T}_2}$, and
$\widetilde{\mathcal{H}}_{\mathcal{T}} = \mathcal{S}'_{\mathcal{T}_1} \cup \mathcal{S}'_{\mathcal{T}_2} \cup \mathcal{S}_{\mathcal{T}_3}$, with cardinalities $t + 1$ and $t + 2$, respectively. The corresponding redundancy is $\rho = (2t + 3)/(2t + 1)$. An illustration of the oversampled bipartite graph construction in the time domain is shown in Fig. \ref{fig03}(a). Compared with conventional critically sampled multi-stage bipartite decompositions, the proposed method enables a single-stage transformation that retains all temporal edge connections. 

\textit{Remark 4:} In bipartite decompositions, temporally adjacent nodes may be disconnected due to partitioning. Adding \emph{vertical edges} between duplicated nodes in the oversampled bipartite graph restores these connections and preserves the temporal and structural integrity of the original graph.

\subsection{Joint Extension of the Time-Vertex Domain}
To support the design of joint filter banks, we propose a novel method for constructing oversampled bipartite graphs in the joint time-vertex domain based on graph coloring strategies. This method independently builds oversampled bipartite structures in both the time and vertex domains, and subsequently combines them via a tensor product. The resulting graph preserves all original adjacency relations while enabling flexible control of redundancy, thereby facilitating efficient joint signal processing.

In the time domain, the temporal structure is modeled as a ring graph, which is always 3-colorable according to graph coloring theory. In the vertex domain, the underlying graph is assumed to be $K$-colorable, as show in Fig. \ref{fig03}(b). In both domains, bipartite oversampled structures are constructed through coloring-based partitioning and node duplication. Building upon these domain-specific extensions, the bipartite structure of the joint graph is obtained through a tensorized composition, whose validity is formally stated as follows.

\textit{Theorem 1:} Let $\widetilde{\mathcal{T}} = (\widetilde{\mathcal{V}}_{\mathcal{T}}, \widetilde{\mathcal{E}}_{\mathcal{T}})$ and $\widetilde{\mathcal{G}} = (\widetilde{\mathcal{V}}_{\mathcal{G}}, \widetilde{\mathcal{E}}_{\mathcal{G}})$ be two bipartite graphs obtained via oversampled coloring-based extensions in the time and vertex domains, respectively, with corresponding bipartitions $\widetilde{\mathcal{V}} = \widetilde{\mathcal{L}}_{\mathcal{G}} \cup \widetilde{\mathcal{H}}_{\mathcal{G}}$, and $\widetilde{\mathcal{V}}_{\mathcal{T}} = \widetilde{\mathcal{L}}_{\mathcal{T}} \cup \widetilde{\mathcal{H}}_{\mathcal{T}}$. Then, the tensor product graph $\widetilde{\mathcal{J}} = \widetilde{\mathcal{T}} \otimes \widetilde{\mathcal{G}}$ admits a valid bipartite structure, where the left and right node sets are given by
\[
\widetilde{\mathcal{L}} = 
(\widetilde{\mathcal{L}}_{\mathcal{T}} \times \widetilde{\mathcal{L}}_{\mathcal{G}}) 
\cup (\widetilde{\mathcal{H}}_{\mathcal{T}} \times \widetilde{\mathcal{H}}_{\mathcal{G}}), 
\]
\[
\widetilde{\mathcal{H}} = 
(\widetilde{\mathcal{L}}_{\mathcal{T}} \times \widetilde{\mathcal{H}}_{\mathcal{G}}) 
\cup (\widetilde{\mathcal{H}}_{\mathcal{T}} \times \widetilde{\mathcal{L}}_{\mathcal{G}}).
\]
\begin{proof}
The proof is reported in Appendix \ref{appendix B}.
\end{proof}

Furthermore, the total vertex set of the joint graph is constructed as
\[
\widetilde{\mathcal{V}}_{\mathcal{T}} \times \widetilde{\mathcal{V}} =
(\widetilde{\mathcal{L}}_{\mathcal{T}} \cup \widetilde{\mathcal{H}}_{\mathcal{T}}) \times 
(\widetilde{\mathcal{L}}_{\mathcal{G}} \cup \widetilde{\mathcal{H}}_{\mathcal{G}}).
\]
This construction guarantees that adjacency information in both the temporal and vertex domains is preserved in the joint graph, enabling consistent spectral representations for joint filtering. The resulting redundancy is multiplicative, given by
\[
\rho = \rho_{\mathcal{T}} \cdot \rho_{\mathcal{G}} =
(2 - \frac{I(\mathcal{L}_{\mathcal{G}}) + I(\mathcal{H}_{\mathcal{G}})}{N})  \cdot  (2 - \frac{I(\mathcal{L}_{\mathcal{T}}) + I(\mathcal{H}_{\mathcal{T}})}{T}).
\]

By construction, $\rho_{\mathcal{T}} \ge 1$ and $\rho_{\mathcal{G}} \ge 1$, implying the lower bound $\rho \ge 1$, with equality attained when both domains are inherently bipartite and no duplication is required. In the time domain, when the ring length is odd, the proposed extension introduces at most a constant number of additional nodes, yielding $\rho_{\mathcal{T}} \to 1$ as $T \to \infty$. Consequently, for large-scale graphs, the overall redundancy remains bounded and asymptotically approaches that of the vertex extension.

It is important to emphasize that the proposed scheme does not seek combinatorial minimality of redundancy. Rather, under the constraint of single-stage bipartite decomposition and full preservation of original adjacency relations, it achieves near-minimal redundancy while enabling structurally consistent joint filter bank design.

\textit{Remark 5:} The proposed joint extension method differs from directly oversampling the joint graph, which treats it as a generic structure. Instead, structural clarity across domains is preserved, enabling theoretically consistent joint filter bank construction.

\section{Experiments and Results}
This section presents experiments on arbitrary graphs, images, and videos, employing the Harary algorithm \cite{Bigraphs,biparticity} to evaluate the performance gains of the joint time-vertex oversampled Laplacian in reconstruction and denoising.

\subsection{Reconstruction of Joint Time-Vertex Graph Signal}

To evaluate the reconstruction performance of the proposed oversampled joint time-vertex filter bank, we conduct epidemic spreading simulations on the Minnesota Road Network \cite{GSPBOX}. The graph consists of $N=2642$ nodes, each representing an intersection, with edges denoting road connections, and is widely adopted for modeling diffusion processes related to mobility and transportation. The epidemic dynamics are modeled using the Susceptible-Exposed-Infectious-Recovered-Susceptible (SEIRS) framework~\cite{JFT,JFRFT,JLCT}, where recovered individuals lose immunity after a finite period. The state variables are denoted as $S(t)$, $E(t)$, $I(t)$, and $R(t)$, representing the numbers of susceptible, exposed, infectious, and temporarily immune individuals, respectively. The discrete-time evolution is governed by
\begin{equation}
	\begin{aligned}
		S(t+1) &= S(t) - \beta S(t) I(t) + \xi R(t), \\
		E(t+1) &= E(t) + \beta S(t) I(t) - \alpha E(t), \\
		I(t+1) &= I(t) + \alpha E(t) - \gamma I(t), \\
		R(t+1) &= R(t) + \gamma I(t) - \xi R(t),
	\end{aligned}
\end{equation}
where $\beta$ denotes the contagion probability, $\alpha=1/T_e$ the incubation rate with latency period $T_e$, $\gamma=1/T_i$ the recovery rate with infectious period $T_i$, and $\xi=1/T_r$ the immunity-loss rate with immunity period $T_r$.

\begin{table*}[h!]
	\centering
	\caption{Reconstruction MSE comparison across different filter banks.}
	\label{tab1}
	\renewcommand{\arraystretch}{1.2} 
	\footnotesize
	\begin{tabular}{c c c c c}
		\hline
		\textbf{Filter Bank} & $\mathbf{X}_1$ & $\mathbf{X}_2$ & $\mathbf{X}_3$ & $\mathbf{X}_4$ \\
		\hline
		\textbf{Graph-QMF} & $4.31 \times 10^{-11}$ & $1.24 \times 10^{-08}$ & $8.82 \times 10^{-12}$ & $1.33 \times 10^{-08}$ \\
		\textbf{GraphBior} & $4.31 \times 10^{-11}$ & $5.25 \times 10^{-09}$ & $5.15 \times 10^{-12}$ & $5.60 \times 10^{-09}$ \\
		\textbf{Critically-sampled} & $3.63 \times 10^{-27}$ & $9.00 \times 10^{-25}$ & $7.07 \times 10^{-28}$ & $9.72 \times 10^{-25}$ \\
		\textbf{Oversampled} & $3.32 \times 10^{-27}$ & $8.80 \times 10^{-25}$ & $6.23 \times 10^{-28}$ & $9.67 \times 10^{-25}$ \\
		\hline
	\end{tabular}
\end{table*}

\begin{figure}[h!]
	\centering
	\includegraphics[width=\linewidth]{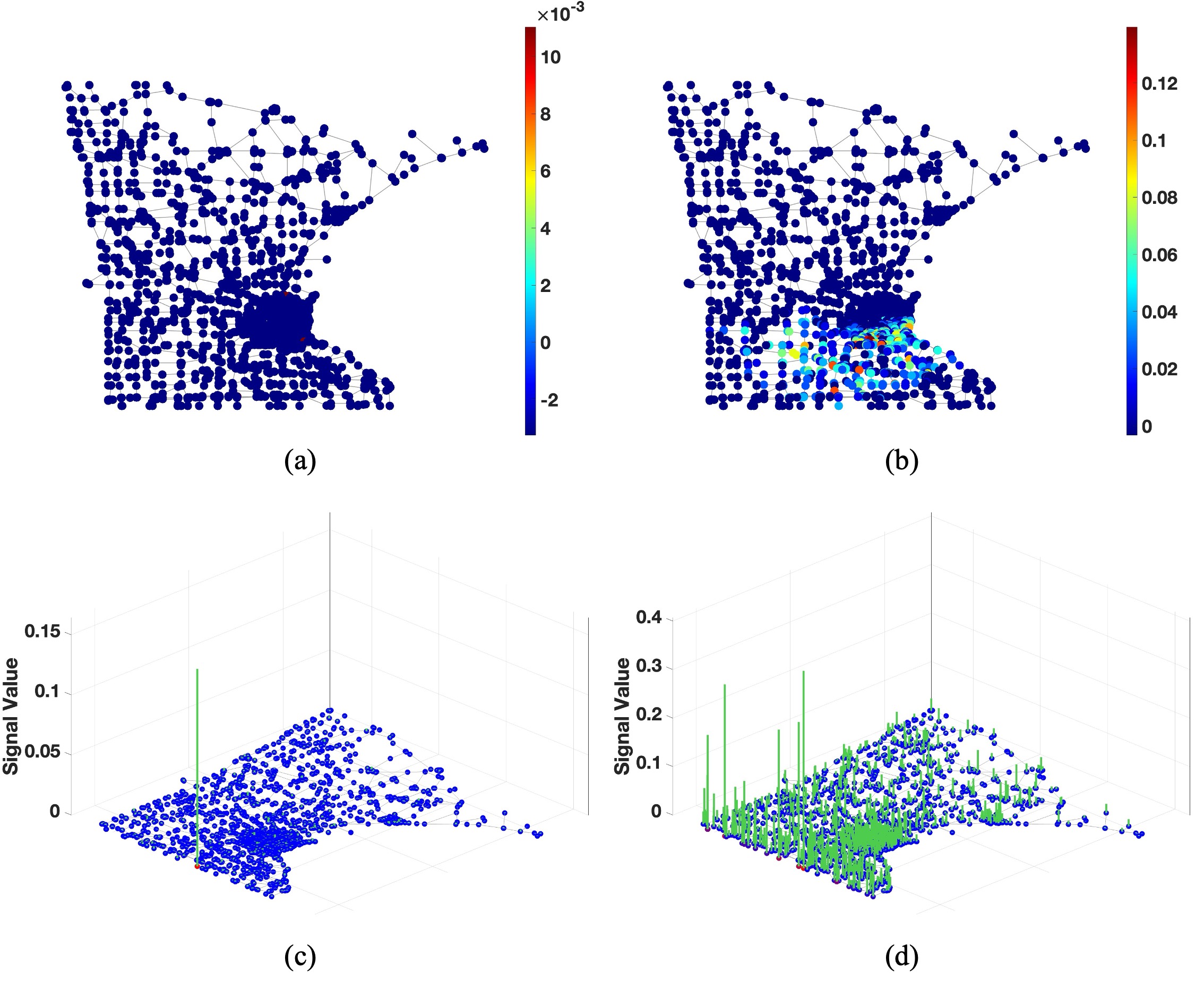}
	\vspace*{-15pt}
	\caption{Joint time-vertex signal $\mathbf{X}_1$ representation of epidemic spreading on the Minnesota Road Network. Panels (a) and (b) show the infection distributions of at the early and late stages, while (c) and (d) depict their corresponding GFT spectrum.}
	\label{fig04}
\end{figure}

Each node hosts 70 individuals, and infections spread across edges with a given contagion probability. Four epidemic scenarios are considered by combining low or high contagion probability with either temporary or permanent immunity. The simulations start with three randomly selected patient-zero nodes, with a latency period of $T_e=2$ days, an infectious period of $T_i=6$ days, and a total duration of approximately 16 months. The normalized infection counts across nodes yield four joint time-vertex signals, denoted as $\{\mathbf{X}_1,\mathbf{X}_2,\mathbf{X}_3,\mathbf{X}_4\}$, corresponding to low contagion with temporary immunity, high contagion with temporary immunity, low contagion with permanent immunity, and high contagion with permanent immunity. Fig.~\ref{fig04} illustrates the spatial distribution of infections for $\mathbf{X}_1$ at the early and late stages of the outbreak, together with their GFT spectrum, which highlight the spectral signatures of the epidemic dynamics.

\begin{figure}[h!]
	\centering
	\includegraphics[width=\linewidth]{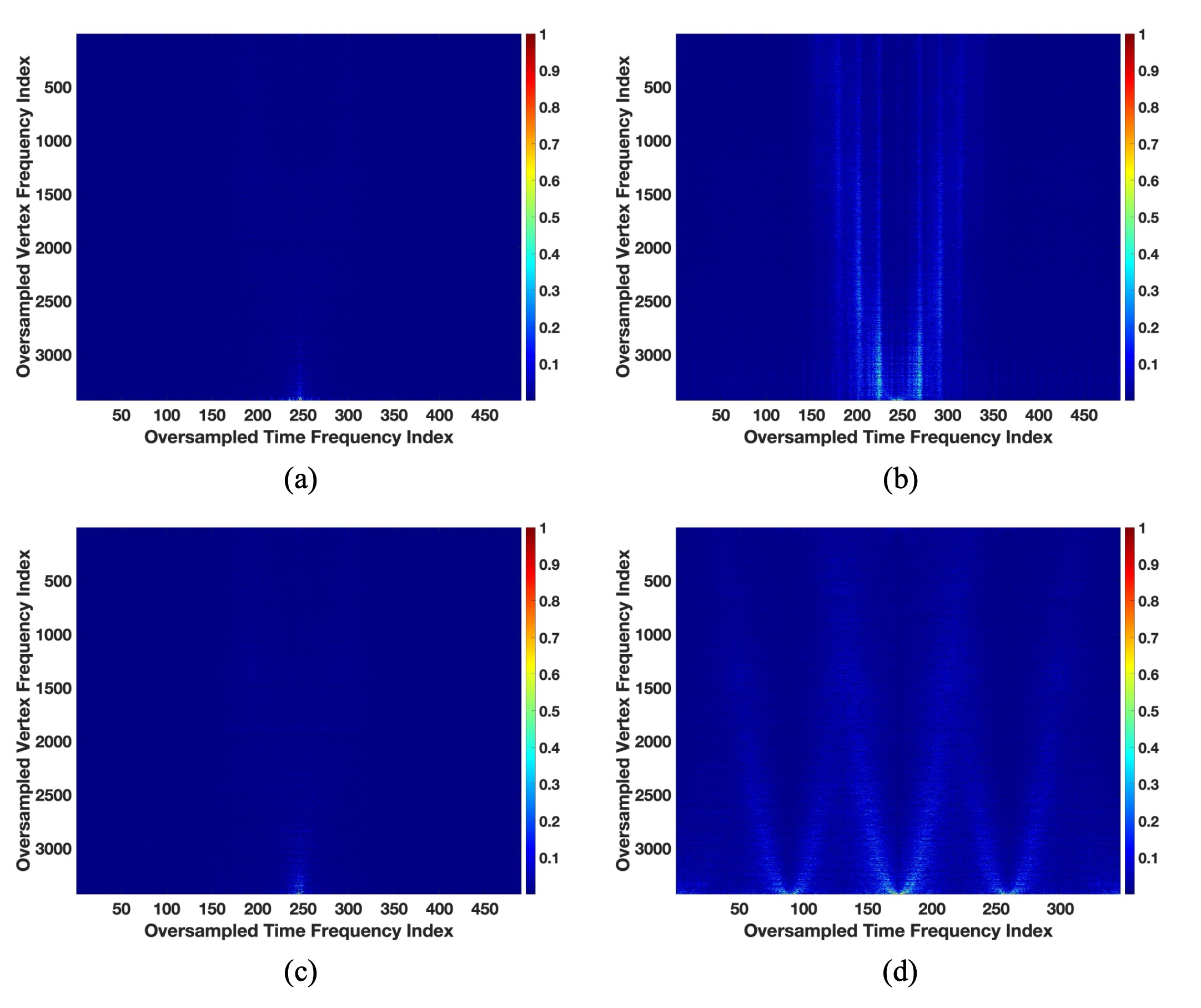}
	\vspace*{-15pt}
	\caption{The JFT of oversampled graph signals of infected populations under different epidemic realizations with varying models and contagion probabilities. (a) Low contagion with temporary immunity. (b) High contagion with temporary immunity. (c) Low contagion with permanent immunity. (d) High contagion with permanent immunity.}
	\label{fig05}
\end{figure}

The joint oversampling procedure is then applied using a $K$-colorable scheme, yielding oversampled graph signals. Their JFT captures the differences in infection dynamics across the four epidemic scenarios, as illustrated in Fig.~\ref{fig05}. 

Finally, reconstruction is performed using the proposed oversampled filter bank and compared against critically-sampled, graph-QMF, and graphBior approaches \cite{biorth,PRwavelet,JFT}. Table~\ref{tab1} reports the reconstruction mean squared error (MSE). The oversampled design consistently achieves significantly lower reconstruction errors, outperforming graph-QMF and graphBior, and critically sampled schemes.

\subsection{Joint Time-Vertex Graph Signal Denoising}
In this section, we evaluate the proposed joint time-vertex oversampled graph filter bank for denoising noisy graph signals and compare it with representative baseline methods. Two distinct types of joint time-vertex models are considered. The first is the Yale Coat of Arms \footnote{[Online]. Available: \url{https://www.cs.yale.edu/homes/spielman/561/}}, which exhibits clear geometric boundaries and local structures. The second is the Minnesota Road Network, which relies on a complex topological structure and better reflects the applicability of the method to sparse graphs. To construct the joint time-vertex signals, the original graph signal is replicated along the temporal dimension to ensure identical input at each time step, followed by the addition of Gaussian noise to generate noisy joint time-vertex signals.

For all methods, the denoising process follows the same procedure. The signal is first decomposed along the temporal dimension, where noise is attenuated by thresholding before reconstruction. The vertex-domain decomposition and thresholding are then applied, and the final reconstruction is obtained in the joint domain. To ensure fairness, identical thresholds are adopted across all methods under different noise levels. Four schemes are tested: the proposed oversampled filter bank, a critically sampled counterpart, the graphBior filter bank, and the graph-QMF filter bank \cite{PRwavelet}. The denoised results are compared against the ground truth, and both MSE and signal-to-noise ratio (SNR) are computed.  

The numerical results are summarized in Table \ref{tab2} and Table \ref{tab3}, corresponding to the Yale Coat of Arms and the Minnesota Road Network, respectively. As shown, the oversampled method consistently achieves the highest SNR across all noise levels, demonstrating superior robustness and fidelity. This improvement stems from the redundancy introduced by oversampling, which enhances the separation between signal and noise in the joint time-vertex domain. The critically-sampled method performs slightly worse but still outperforms the graph-QMF and graphBior approaches.
\begin{table}[h!]
	\centering
	\caption{Yale Coat of Arms denoising results: SNR (dB)}
	\label{tab2}
	\renewcommand{\arraystretch}{1.2}
	\footnotesize
	\begin{tabular}{c ||c |c |c |c ||c}
		\hline
		$\sigma$ & 1/8 & 1/4 & 1/2 & 1 & $\rho$ \\
		\hline
		\textbf{Noisy} & 17.90 & 11.88 & 5.86 & -0.16& -- \\
		\textbf{Graph-QMF} & 15.23 & $12.39$ & $9.04$ & $5.02$& 1.00 \\
		\textbf{GraphBior} & 14.40 & $12.31$ & $8.44$ & $4.48$ & 1.00 \\
		\textbf{Critically sampled} & 14.49 & $12.68$ & $9.27$ & $5.71$ & 1.00 \\
		\textbf{Oversampled} & \textbf{19.92} & $\textbf{15.83}$ & $\textbf{11.16}$ & $\textbf{6.36}$ & 3.05\\
		\hline
	\end{tabular}
\end{table}

\begin{table}[h!]
	\centering
	\caption{Minnesota Road Network denoising results: SNR (dB)}
	\label{tab3}
	\renewcommand{\arraystretch}{1.2}
	\footnotesize
	\begin{tabular}{c ||c |c |c |c ||c}
		\hline
		$\sigma$ & 1/8 & 1/4 & 1/2 & 1 & $\rho$ \\
		\hline
		\textbf{Noisy} & 18.03 & 12.01 & 5.99 & -0.03& -- \\
		\textbf{Graph-QMF} & 20.39 & 12.69 & 9.62 & 5.29& 1.00 \\
		\textbf{GraphBior} & 19.52 & 13.56 & 8.29 & 3.01 & 1.00 \\
		\textbf{Critically sampled} & 20.75 & 13.72 & 10.39 & 5.51 & 1.00 \\
		\textbf{Oversampled} & \textbf{20.91} & \textbf{14.52} & \textbf{10.67} & \textbf{5.53} & 1.30\\
		\hline
	\end{tabular}
\end{table}

\begin{figure*}[h!]
	\centering
	\includegraphics[scale=0.15]{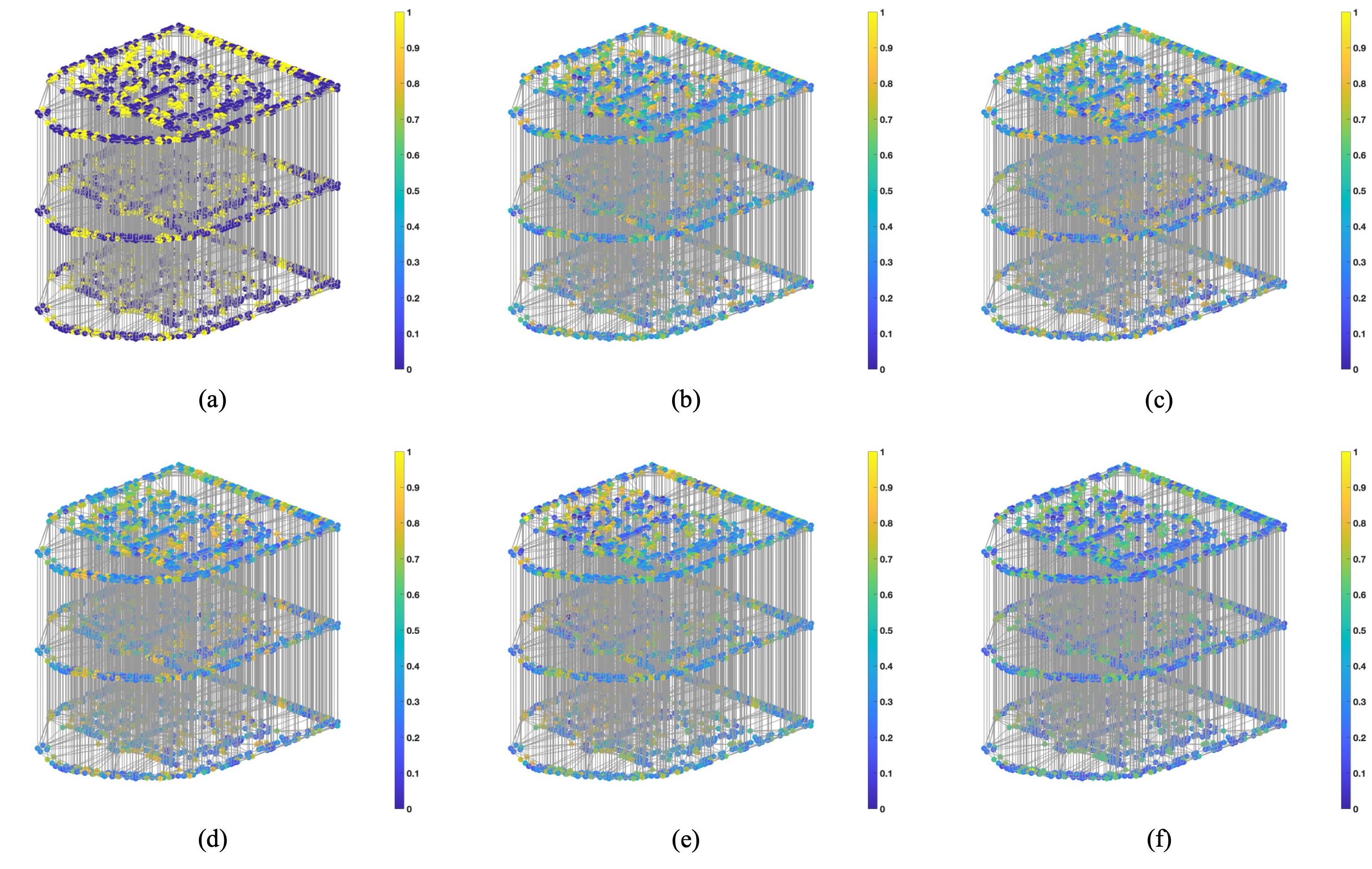}
	\vspace*{-3pt}
	\caption{Denoising results of the Yale Coat of Arms: (a) original signal; (b) noisy signal; (c) Graph-QMF; (d) GraphBior; (e) critically sampled filter bank; (f) oversampled filter bank.}
	\label{fig06}
\end{figure*}

\begin{figure*}[h!]
	\centering
	\includegraphics[scale=0.15]{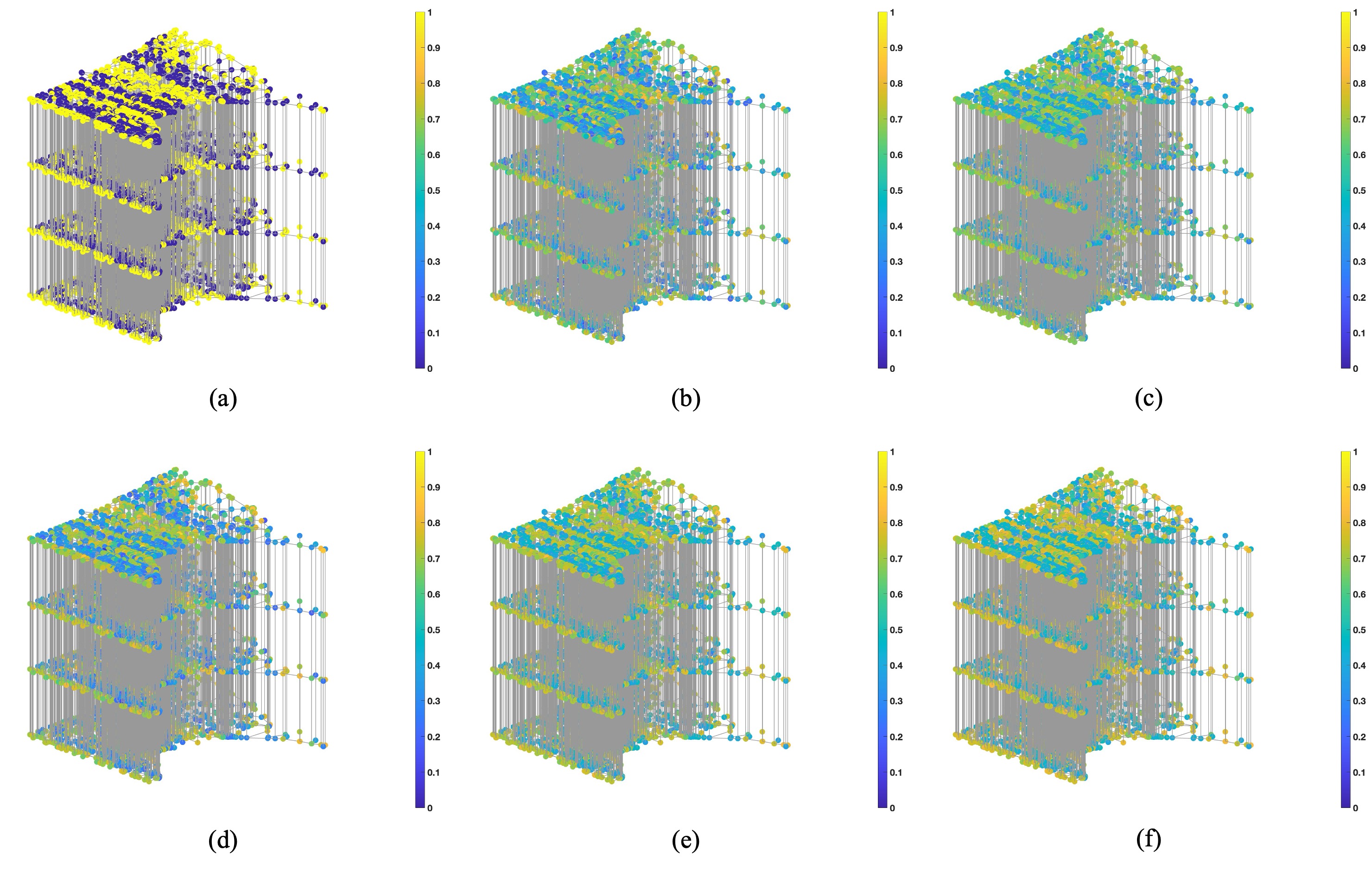}
	\vspace*{-3pt}
	\caption{Denoising results of the Minnesota Road Network: (a) original signal; (b) noisy signal; (c) Graph-QMF; (d) GraphBior; (e) critically-sampled filter bank; (f) oversampled filter bank.}
	\label{fig07}
\end{figure*}

\begin{figure}[h!]
	\centering
	\includegraphics[width=\linewidth]{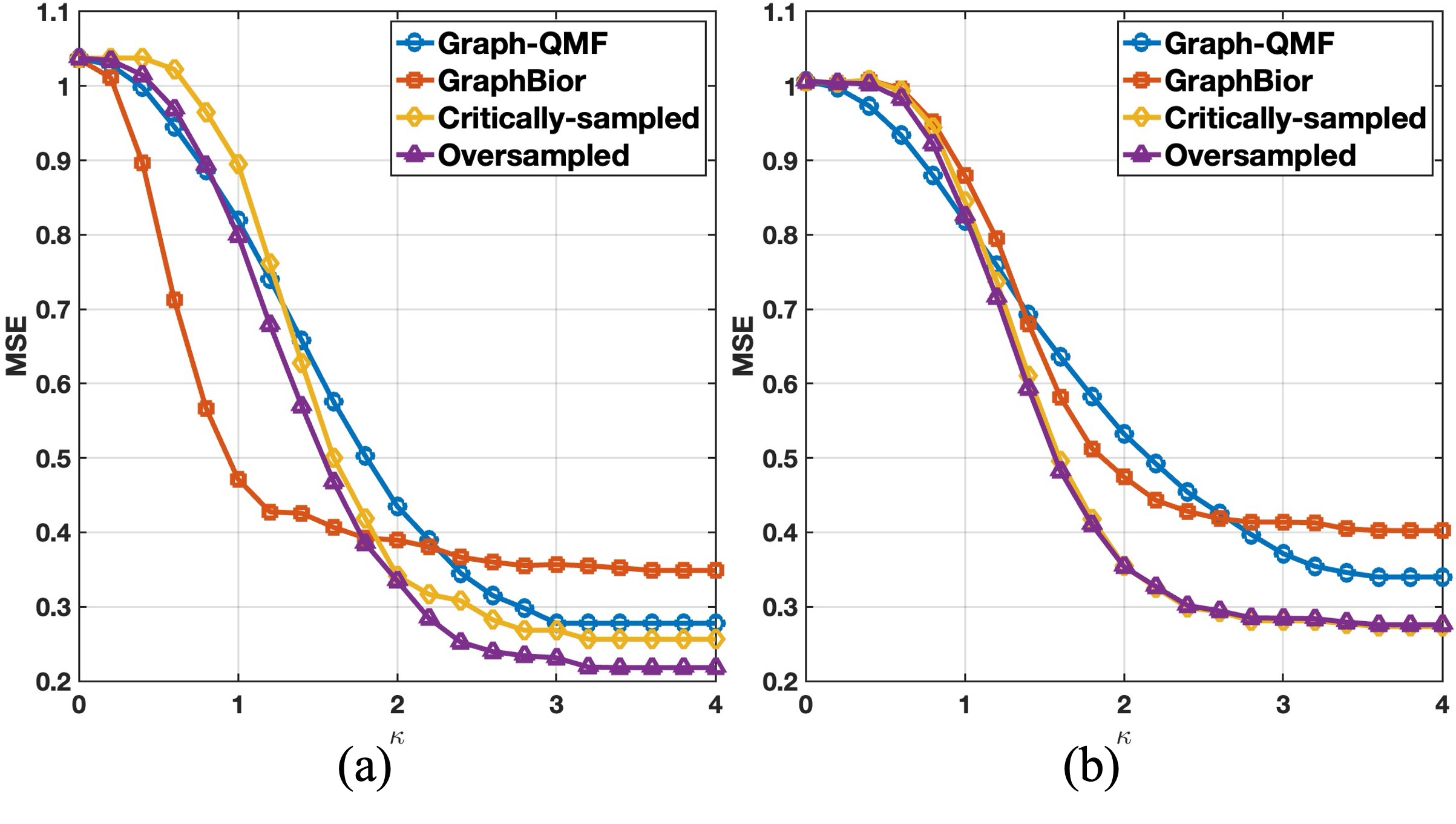}
	\vspace*{-15pt}
	\caption{MSE as a Function of the Threshold Scaling Coefficient $\kappa$: (a) Yale Coat of Arms; (b) Minnesota Road Network.}
	\label{fig08}
\end{figure}

To provide further insights, Fig.~\ref{fig06} and Fig.~\ref{fig07} present qualitative comparisons of the denoising results for the Yale Coat of Arms and the Minnesota Road Network under Gaussian noise with $\sigma = 1/2$, respectively.
Fig.~\ref{fig08} further reports the quantitative comparison of the denoising performance in terms of MSE for (a) the Yale Coat of Arms and (b) the Minnesota Road Network under Gaussian noise with $\sigma = 1$. The results are evaluated over a range of threshold values defined as $\mathrm{threshold} = \kappa \sigma$, where $\kappa$ denotes the scaling coefficient. It can be observed that, when the threshold parameter is properly tuned, the proposed oversampled scheme consistently achieves the lowest reconstruction error among the compared methods. Overall, the oversampled filter bank exhibits clear advantages in joint time-vertex signal denoising, highlighting its strong potential for complex and large-scale signal processing applications.

\subsection{Image Denoising via Joint Oversampled Filter Banks}
To evaluate the efficacy of the proposed two-channel joint time-vertex oversampled graph filter banks, we conduct experiments on both static images and dynamic videos. Within the joint time-vertex framework, an image or a sequence of frames is modeled as a graph, where pixels are treated as vertices and temporal evolution is captured along the time dimension. This allows simultaneous filtering across spatial and temporal domains.

For static image denoising, two pairs of handwritten digit images are selected: \textit{0–6} and \textit{1–7}. Each pair serves as the endpoints of a joint sequence, and $T=101$ intermediate frames are generated by linear interpolation along the temporal axis, forming a data matrix. Each image is resized to $N \times N$ ($N=35$ in this study), with columns representing temporal frames and rows corresponding to graph vertices (pixels). Gaussian noise with standard deviation $\sigma=20$ is added to construct noisy image sequences. Both oversampled and critically sampled filter banks are applied for denoising, and the first and last frames are visualized for comparison.

\begin{figure}[h]
	\begin{center}
		\begin{minipage}[t]{0.9\linewidth}
			\centering
			\includegraphics[width=\linewidth]{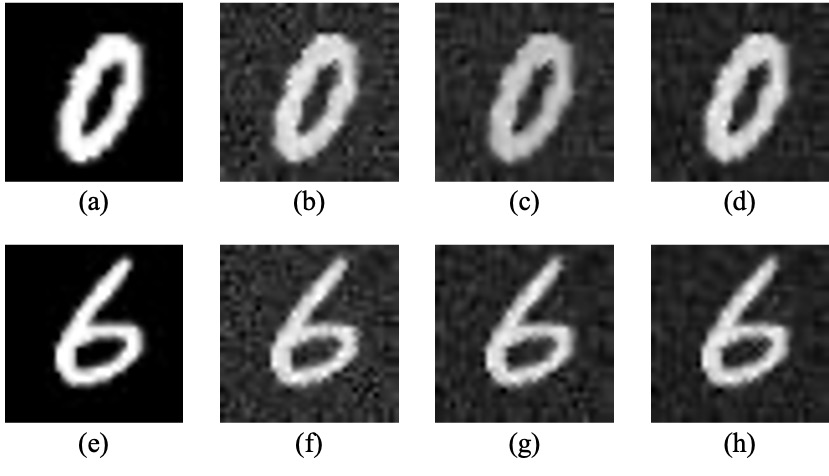}
		\end{minipage}
	\end{center}
	\caption{Image denoising results for digit pair \textit{0--6}: (a) original 0, (b) noisy 0, (c) critically sampled reconstruction, (d) oversampled reconstruction, (e) original 6, (f) noisy 6, (g) critically sampled reconstruction, (h) oversampled reconstruction.}
	\vspace*{-3pt}
	\label{fig09}
\end{figure}

\begin{figure}[h]
	\begin{center}
		\begin{minipage}[t]{0.9\linewidth}
			\centering
			\includegraphics[width=\linewidth]{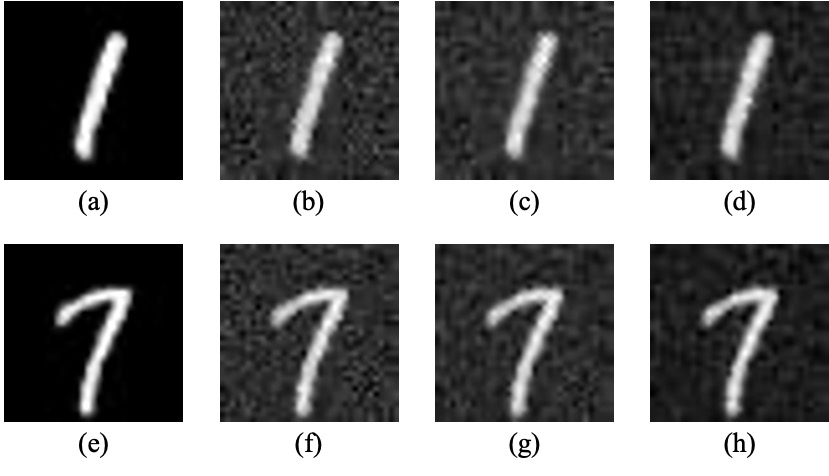}
		\end{minipage}
	\end{center}
	\caption{Image denoising results for digit pair \textit{1--7}: 
		(a) original 1, (b) noisy 1, (c) critically sampled reconstruction, (d) oversampled reconstruction, (e) original 7, (f) noisy 7, (g) critically sampled reconstruction, (h) oversampled reconstruction.}
	\vspace*{-3pt}
	\label{fig10}
\end{figure}

Results are shown in Figs.~\ref{fig09} and~\ref{fig10}, where panels (a) and (e) correspond to the original images, (b) and (f) to the noisy inputs, (c) and (g) to reconstructions using critically sampled filter banks, and (d) and (h) to reconstructions using oversampled filter banks. Visually, the oversampled method more effectively suppresses noise while preserving edges and fine details. Quantitative evaluations under different noise levels are summarized in Tables~\ref{tab4} and~\ref{tab5}, which report SNR improvements for both image pairs.

\begin{table}[h!]
	\centering
	\caption{Digit pair \textit{0--6} denoising: SNR (dB)}
	\label{tab4}
	\renewcommand{\arraystretch}{1.2}
	\footnotesize
	\begin{tabular}{c ||c |c |c |c ||c}
		\hline
		$\sigma$ & 1 & 5 & 10 & 20 & $\rho$ \\
		\hline
		\textbf{Noisy} & 40.54 & 26.56 & 20.54 & 14.52 & -- \\
		\textbf{Critically sampled} & 40.72 & 26.80 & 21.25 & 16.27 & 1.00 \\
		\textbf{Oversampled} & \textbf{40.79} & \textbf{26.90} & \textbf{21.42} & \textbf{16.62} & 1.02 \\
		\hline
	\end{tabular}
\end{table}

\begin{table}[h!]
	\centering
	\caption{Digit pair \textit{1--7} denoising: SNR (dB)}
	\label{tab5}
	\renewcommand{\arraystretch}{1.2}
	\footnotesize
	\begin{tabular}{c ||c |c |c |c ||c}
		\hline
		$\sigma$ & 1 & 5 & 10 & 20 & $\rho$ \\
		\hline
		\textbf{Noisy} & 36.66 & 22.68 & 16.66 & 10.64 & -- \\
		\textbf{Critically sampled} & 37.32 & 23.55 & 17.94 & 12.88 & 1.00 \\
		\textbf{Oversampled} & \textbf{37.36} & \textbf{23.62} & \textbf{18.11} & \textbf{13.15} & 1.02 \\
		\hline
	\end{tabular}
\end{table}

\begin{figure*}[h!]
	\centering
	\includegraphics[scale=0.15]{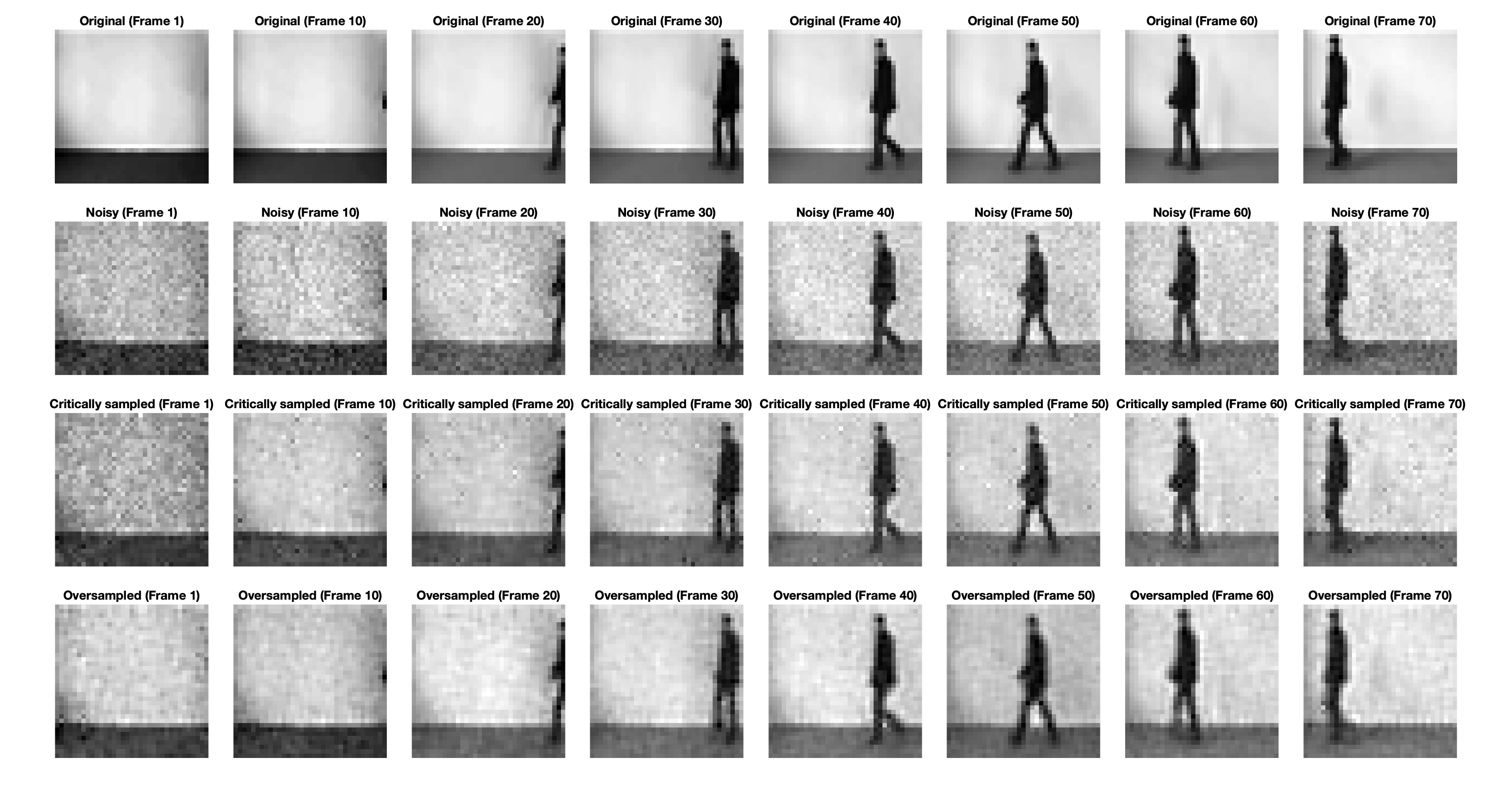}
	\vspace*{-3pt}
	\caption{Key-frame denoising results for the walking sequence: original frames (row 1), noisy frames with $\sigma=10$ (row 2), critically sampled reconstructions (row 3), and oversampled reconstructions (row 4). Frames shown: 1, 10, 20, 30, 40, 50, 60, 70.}
	\label{fig11}
\end{figure*}

\begin{figure}[h!]
		\centering
		\includegraphics[width=\linewidth]{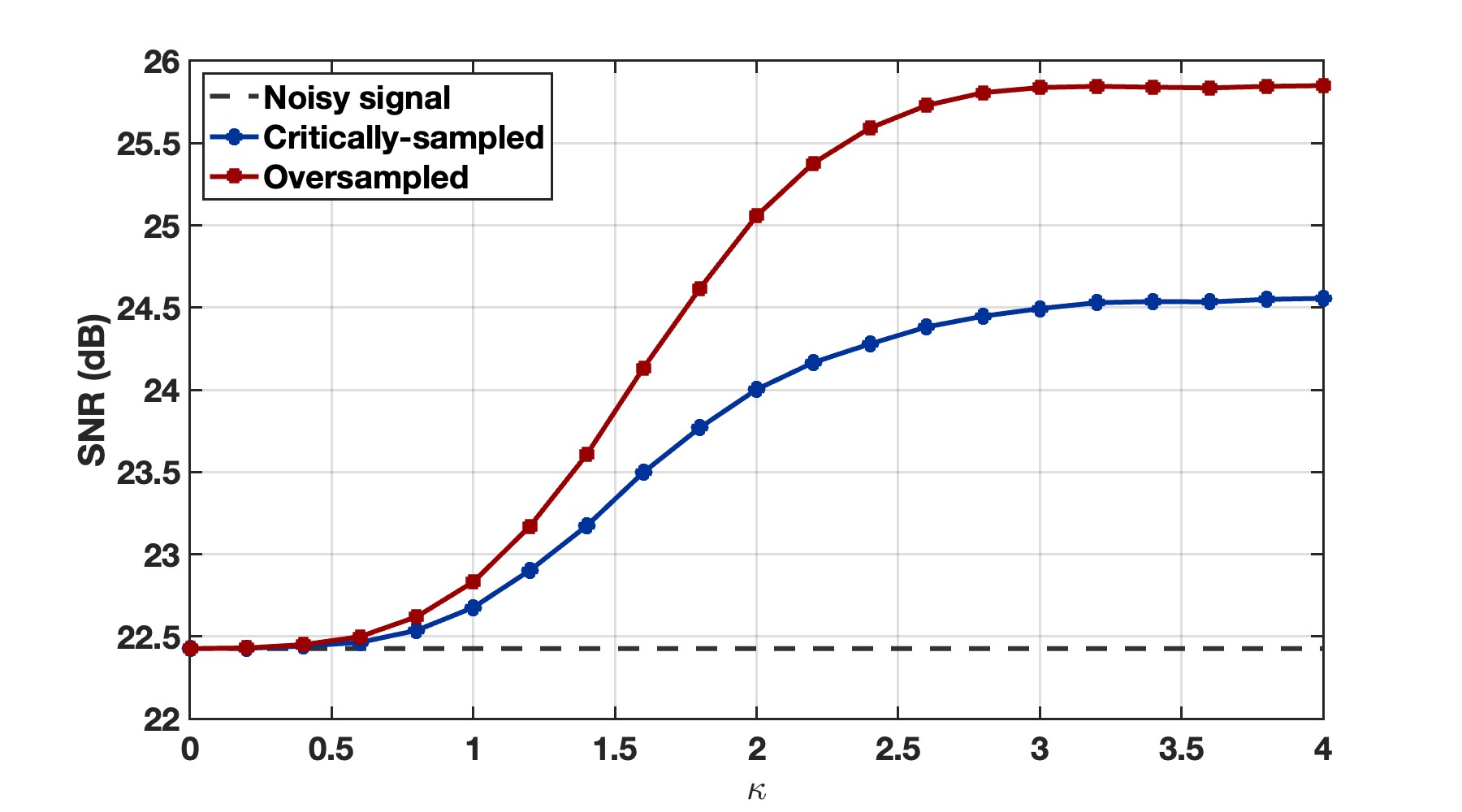}
		\vspace*{-3pt}
		\caption{SNR as a function of the threshold scaling coefficient $\kappa$ for video denoising.}
		\label{fig12}
\end{figure}

For dynamic video denoising, the walking sequence from the human actions video database\footnote{[Online]. Available: \url{https://www.csc.kth.se/cvap/actions/}} is used. The first 3 seconds (approximately 75 frames) are extracted, converted to grayscale, and resized to $N \times N$, forming a spatial graph with each pixel as a node. The temporal dimension is represented by a ring graph over the frame sequence, yielding a joint time-vertex graph structure. Gaussian noise with standard deviation $\sigma=10$ is added to the frames, and both oversampled and critically sampled filter banks are applied for denoising. Selected key frames are plotted in Fig.~\ref{fig11}, with rows corresponding to original frames, noisy frames, critically sampled reconstructions, and oversampled reconstructions. Fig.~\ref{fig12} illustrates the influence of the threshold scaling coefficient $\kappa$ on the video denoising performance. It can be observed that the proposed oversampled scheme consistently achieves the highest SNR across the entire range of $\kappa$ values.
The noisy video yields an SNR of 22.43 dB, which increases to 24.55 dB with the critically sampled filter bank and further to 25.85 dB using the oversampled filter bank. The corresponding execution times are 0.79 s for the critically sampled scheme and 1.07 s for the oversampled scheme. The performance gain stems from the redundancy introduced by the oversampled filter bank, which enhances reconstruction accuracy and denoising capability while better preserving motion details and temporal continuity. 

This improvement comes at the expense of higher computational cost compared to the critically sampled design. The runtime increases proportionally with the oversampling factor $\rho$, as well as with the graph size and temporal dimension. Nevertheless, for moderate problem scales, the computational overhead remains acceptable and is justified by the significant performance gain.

\section{Conclusion}
In this work, we proposed a novel two-channel filter bank for joint time-vertex signal processing. An oversampled graph Laplacian operator was formulated and extended from vertex-domain filter banks to a joint oversampled graph filter bank. Based on this operator, a two-channel oversampled filter bank was designed using a $K$-coloring strategy, enabling the construction of oversampled bipartite graphs that preserve all temporal and spatial edges. Extensive experiments on both signal reconstruction and image denoising validate the effectiveness of the proposed filter banks, demonstrating notable improvements in denoising performance compared to existing methods. Future work will focus on optimized and adaptive filter design, extensions to directed and dynamic graphs, and broader applications in large-scale time-varying signal data.
	
		
\appendices
\section{Colorable Graph Extension Methods}
\label{appendix A}
To construct oversampled bipartite graphs suitable for filter bank design, a general method based on graph coloring is adopted \cite{OGFB,OGLM,OGFT}. Consider a $K$-colorable graph $\mathcal{G} = (\mathcal{V}, \mathcal{E})$ with $K \geq 3$. The objective is to generate a bipartite extension while preserving all original edges.

Assume the vertex set is partitioned into disjoint color classes
\[
\mathcal{V} = \mathcal{S}_{\mathcal{G}_1} \cup \mathcal{S}_{\mathcal{G}_2} \cup \dots \cup \mathcal{S}_{\mathcal{G}_K},
\]
such that adjacent nodes have different colors. For an integer $1 \leq l < K$, define the left and right node sets of the foundation bipartite graph $\mathcal{G}_b$ as
\begin{equation}
	\mathcal{L}_{\mathcal{G}_b} = \bigcup_{i=1}^l \mathcal{S}_{\mathcal{G}_i}, \quad \mathcal{H}_{\mathcal{G}_b} = \bigcup_{i=l+1}^K \mathcal{S}_{\mathcal{G}_i},
\end{equation}
and retain only edges $\mathcal{E}_b$ between $\mathcal{L}_{\mathcal{G}_b}$ and $\mathcal{H}_{\mathcal{G}_b}$.

To preserve the remaining intra-set edges while maintaining bipartiteness, duplicated nodes are introduced as follows:
\begin{itemize}
	\item For each $\mathcal{S}_{\mathcal{G}_i} \subseteq \mathcal{L}_{\mathcal{G}_b}$, a duplicate $\mathcal{S}'_{\mathcal{G}_i}$ is added to the opposite set $\widetilde{\mathcal{H}}_{\mathcal{G}}$.
	\item For each $\mathcal{S}_{\mathcal{G}_j} \subseteq \mathcal{H}_{\mathcal{G}_b}$, a duplicate $\mathcal{S}'_{\mathcal{G}_j}$ is added to $\widetilde{\mathcal{L}}_{\mathcal{G}}$.
\end{itemize}

Each duplicated node is vertically connected to its original, and all residual edges are reconstructed between bipartite layers. The resulting oversampled bipartite graph $\widetilde{\mathcal{G}}$ has node sets 
\[
\widetilde{\mathcal{L}}_{\mathcal{G}} = \bigcup_{i=1}^l \mathcal{S}_{\mathcal{G}_i} \cup \bigcup_{j=l+1}^K \mathcal{S}'_{\mathcal{G}_j}, \quad \widetilde{\mathcal{H}}_{\mathcal{G}} = \bigcup_{i=1}^l \mathcal{S}'_{\mathcal{G}_i}  \cup \bigcup_{j=l+1}^K \mathcal{S}_{\mathcal{G}_j}.
\]

The number of isolated nodes in $\mathcal{L}_{\mathcal{G}_b}$ and $\mathcal{H}_{\mathcal{G}_b}$, denoted $I(\mathcal{L}_{\mathcal{G}_b})$ and $I(\mathcal{H}_{\mathcal{G}_b})$ respectively, do not require duplication. The total sizes of $\widetilde{\mathcal{L}}_{\mathcal{G}} $ and $\widetilde{\mathcal{G}}_{\mathcal{G}} $ are
\[
|\widetilde{\mathcal{L}}_{\mathcal{G}}| = N - I(\mathcal{H}_{\mathcal{G}_b}), \quad
|\widetilde{\mathcal{H}}_{\mathcal{G}}| = N - I(\mathcal{L}_{\mathcal{G}_b}),
\]
where $N = |\mathcal{V}|$. The redundancy for the two-channel oversampled filter bank is given by
\[
\rho = 2 - \frac{I(\mathcal{L}_{\mathcal{G}_b})+I(\mathcal{H}_{\mathcal{G}_b})}{N}.
\]

For three-colorable construction is a special case under $K = 3$ and $l = 2$, yielding
$\mathcal{L}_{\mathcal{G}_b} = \mathcal{S}_{\mathcal{G}_1} \cup \mathcal{S}_{\mathcal{G}_2}$, and $\mathcal{H}_{\mathcal{G}_b} = \mathcal{S}_{\mathcal{G}_3}$. Thus, 
\begin{equation}
	\widetilde{\mathcal{L}}_{\mathcal{G}} = \mathcal{S}_{\mathcal{G}_1} \cup \mathcal{S}_{\mathcal{G}_2}, \quad \widetilde{\mathcal{H}}_{\mathcal{G}} = \mathcal{S}'_{\mathcal{G}_1} \cup \mathcal{S}'_{\mathcal{G}_2} \cup \mathcal{S}_{\mathcal{G}_3}.
\end{equation}

For the bipartite double cover \cite{OGLM,Bigraphs}, defined as the tensor product $\widetilde{\mathcal{G}}_{BDC} = \mathcal{G} \otimes K_2$, with the complete graph of two vertices $K_2$, generates a bipartite graph with $2N$ nodes and $2|\mathcal{E}|$ edges. This corresponds to the case $l = K$ without vertical edges, producing a spectrum with two symmetric copies of the original, resulting in full redundancy. 

In contrast, the color-based method reduces redundancy and enables spectral control, offering greater flexibility for oversampled graph filter bank design.

		
\section{Proof of Theorem 1}\label{appendix B}
By definition of the tensor product graph, the vertex set of the joint graph $\widetilde{\mathcal{T}} \otimes \widetilde{\mathcal{G}}$ is given by the Cartesian product, 
$\widetilde{\mathcal{V}}_{\mathcal{T}} \times \widetilde{\mathcal{V}}$.

Two nodes $(t_1, v_1)$ and $(t_2, v_2)$ are connected in $\widetilde{\mathcal{J}}$ if and only if: $t_1 \sim t_2$ in $\widetilde{\mathcal{T}}$, and $v_1 \sim v_2$ in $\widetilde{\mathcal{G}}$.

Since both $\widetilde{\mathcal{T}}$ and $\widetilde{\mathcal{G}}$ are bipartite, their edges only connect nodes from opposite partitions
$$
\widetilde{\mathcal{V}}_{\mathcal{T}} = \widetilde{\mathcal{L}}_{\mathcal{T}} \cup \widetilde{\mathcal{H}}_{\mathcal{T}}, \quad
\widetilde{\mathcal{V}} = \widetilde{\mathcal{L}}_{\mathcal{G}} \cup \widetilde{\mathcal{H}}_{\mathcal{G}}.
$$

To ensure that every edge in $\widetilde{\mathcal{J}}$ connects nodes from opposite classes, we define the joint bipartition as
$$
\widetilde{\mathcal{L}} := (\widetilde{\mathcal{L}}_{\mathcal{T}} \times \widetilde{\mathcal{L}}_{\mathcal{G}}) 
\cup (\widetilde{\mathcal{H}}_{\mathcal{T}} \times \widetilde{\mathcal{H}}_{\mathcal{G}}),
$$
$$
\widetilde{\mathcal{H}} := (\widetilde{\mathcal{L}}_{\mathcal{T}} \times \widetilde{\mathcal{H}}_{\mathcal{G}}) 
\cup (\widetilde{\mathcal{H}}_{\mathcal{T}} \times \widetilde{\mathcal{L}}_{\mathcal{G}}).
$$

Under this assignment, all adjacent node pairs $(t_1, v_1) \sim (t_2, v_2)$ lie across the two sets $\widetilde{\mathcal{L}}$ and $\widetilde{\mathcal{H}}$, satisfying the bipartite property. Hence, $\widetilde{\mathcal{J}}$ is a valid bipartite graph.

\section*{Acknowledgments}
In the spirit of reproducible research, MATLAB code examples are available at https://github.com/Zhangyubit/JOGLM.

		\newpage
		
		
		
		
		

		
	\end{document}